\documentclass{elsarticle}
\usepackage{amsmath,amssymb,amsthm}

\begin{document}

\begin{frontmatter}

\title{A Framework for Generalising the Newton Method and Other
Iterative Methods from Euclidean Space to Manifolds}

\author{Jonathan H. Manton}
\address{Department of Electrical and Electronic Engineering \\
The University of Melbourne, Victoria 3010, Australia}
\ead{j.manton@ieee.org}

\begin{abstract}
The Newton iteration is a popular method for minimising a cost
function on Euclidean space.  Various generalisations to cost
functions defined on manifolds appear in the literature.  In each
case, the convergence rate of the generalised Newton iteration
needed establishing from first principles.  The present paper
presents a framework for generalising iterative methods from Euclidean
space to manifolds that ensures local convergence rates are preserved.
It applies to any (memoryless) iterative method computing a coordinate
independent property of a function (such as a zero or a local
minimum).  All possible Newton methods on manifolds are believed
to come under this framework.  Changes of coordinates, and not any
Riemannian structure, are shown to play a natural role in lifting
the Newton method to a manifold.  The framework also gives new
insight into the design of Newton methods in general.
\end{abstract}


\journal{Numerische Mathematik (Accepted on 4 April 2014)}

\end{frontmatter}

\newcommand{\Ct}{\mathcal{C}}
\newcommand{\reals}{\mathbb{R}}
\newcommand{\Nplus}{\mathbb{N}_+}	
\newcommand{\rn}{\reals^n}
\newcommand{\grad}{\operatorname{grad}}
\newcommand{\Exp}{\operatorname{Exp}}
\newcommand{\Tr}{\operatorname{Tr}}
\newcommand{\cmp}[1]{{#1}^{\mathrm{c}}}	
\newcommand{\bd}{\partial}		

\newcommand{\Hess}{H}		

\newcommand{\xstar}{{x^{\ast}}}
\newcommand{\xstarf}[1]{{x_{#1}^{\ast}}}
\newcommand{\ystar}{{y^{\ast}}}
\newcommand{\pstar}{{p^{\ast}}}

\newcommand{\Hi}[1]{H1${}_{#1}$}
\newcommand{\Hii}[1]{H2${}_{#1}$}
\newcommand{\Hiii}[1]{H3${}_{#1}$}

\newcommand{\refeq}[1]{(\ref{eq:#1})}
\newcommand{\refsec}[1]{Section~\ref{sec:#1}}
\newcommand{\refapp}[1]{\ref{app:#1}}
\newcommand{\refdef}[1]{Definition~\ref{def:#1}}
\newcommand{\reflem}[1]{Lemma~\ref{lem:#1}}
\newcommand{\refrem}[1]{Remark~\ref{rem:#1}}
\newcommand{\refprop}[1]{Proposition~\ref{prop:#1}}
\newcommand{\refth}[1]{Theorem~\ref{th:#1}}
\newcommand{\refex}[1]{Example~\ref{ex:#1}}
\newcommand{\refcor}[1]{Corollary~\ref{cor:#1}}


\theoremstyle{plain}
\newtheorem{thm}{Theorem}
\newtheorem{cor}[thm]{Corollary}
\newtheorem{prop}[thm]{Proposition}
\newtheorem{lem}[thm]{Lemma}

\theoremstyle{definition}
\newtheorem{defn}[thm]{Definition}
\newtheorem{conjecture}[thm]{Conjecture}
\newtheorem{exmp}[thm]{Example}

\theoremstyle{remark}
\newtheorem{rem}[thm]{Remark}


\section{Introduction}

The Newton iteration function $N_f \colon \rn \rightarrow \rn$
associated with a smooth cost function $f \colon \rn \rightarrow
\reals$ is
\begin{equation}
\label{eq:newt}
N_f(x) = x - [\Hess_f(x)]^{-1} \nabla f(x),\qquad x \in \rn
\end{equation}
where $\nabla f(x)$ and $\Hess_f(x)$ are the gradient and Hessian
of $f$, respectively; $N_f$ does not depend on the choice of inner
product with respect to which the gradient and Hessian are defined.
Starting with an initial guess $x_0 \in \rn$, the Newton method
uses the Newton iteration function to generate the iterates $x_{k+1}
= N_f( x_k )$.  Under certain conditions~\cite{bk:Polak:opt}, this
sequence is well-defined and converges to a critical point of $f$,
meaning $[\Hess_f(x_k)]^{-1}$ exists for all $k$, and $x = \lim_{k
\rightarrow \infty} x_k$ exists and satisfies $\nabla f(x) = 0$.

Let $f\colon M \rightarrow \reals$ now be a smooth cost function defined
on an $n$-dimensional manifold $M$.  Since $M$ locally looks like
$\rn$, it is natural to ask how the Newton iteration function
\refeq{newt} can be extended to an iteration function $E_f\colon M
\rightarrow M$ such that the iterates $x_{k+1} = E_f( x_k )$ enjoy
the same locally quadratic rate of convergence as do the Euclidean
Newton iterates.

One approach~\cite{Gabay:optm} is to endow the manifold $M$ with a
metric and define $E_f$ by a formula analogous to \refeq{newt} but
with $\nabla f$ and $\Hess_f$ replaced by the Riemannian gradient
$\grad f$ and Hessian $\mathcal{H}_f$ of $f$, and the straight-line
increment $-[\Hess_f(x)]^{-1} \nabla f(x)$ replaced by an increment
along a geodesic, namely
\begin{equation}
\label{eq:riem_n}
E_f(p) = \Exp_p\left( -[\mathcal{H}_f(p)]^{-1} \grad f(p) \right)
\end{equation}
where $\Exp_p$ is the Riemannian exponential map centred at $p$.

The Riemannian Newton method \refeq{riem_n} has some disadvantages
and other Newton methods on manifolds are possible~\cite{Adler:2002cc,
Manton:opt_mfold}.

\emph{What is the most general form of a Newton method on a manifold?}
Here, a Newton method is defined as any iterative algorithm $p_{k+1}
= E_f(p_k)$ that converges locally quadratically to every non-degenerate
critical point of every reasonable cost function $f$, where $E_f(p)$
depends only on the 2-jet of the function $f$ at $p$; if $f$ and
$g$ agree to second order at $p$ then $E_f(p) = E_g(p)$.

\refth{nmm} affords an answer, expressed in terms of parametrisations.
A parametrisation of a manifold $M$ is a function $\phi\colon TM
\rightarrow M$ whose restriction $\phi_p\colon T_pM \rightarrow M$
to the tangent space $T_pM$ of $M$ at any point $p \in M$ provides
a (not necessarily one-to-one) correspondence between a neighbourhood
of $0_p \in T_pM$ and a neighbourhood of $p \in M$; the former is
a subset of a vector space and therefore easier to work with.  It
suffices for $\phi$ to be $\Ct^2$-smooth and satisfy $\phi_p(0_p)=p$,
but interestingly, there exist valid parametrisations that are not
continuous.  (Precise definitions are given in the body of the paper.)

\refth{nmm} states that for any pair of parametrisations $\phi\colon
TM \rightarrow M$ and $\psi\colon TM \rightarrow M$, the iteration
function $E_f(p) = \psi_p \circ N_{f \circ \phi_p}(0_p)$ is a Newton
method on the manifold $M$, where $N$ is the Euclidean Newton
iteration function \refeq{newt} but on the abstract vector space
$T_pM$ rather than $\reals^n$.  (Since \refeq{newt} does not depend
on the choice of inner product, there is no need for a Riemannian
metric on $M$.) Justification is given in the body of the paper for
believing this to be the most general form possible of a Newton
method on a manifold.

Requiring a Newton method to be strictly of the form $p_{k+1} =
E_f(p_k)$ places an unnecessary global topological constraint on
the parametrisations.  Instead, $\phi_p\colon T_pM \rightarrow M$
could be constructed on demand by ``transporting'' the old
parametrisation $\phi_{p_{k-1}}$ from $p_{k-1}$ to $p_k$.  As
transport is generally path dependent, $E_f(p_k)$ may
depend on where $p_k$ is relative to $p_{k-1}$.  A uniformity
constraint on the family of possible parametrisations allows for
the generalisation of \refth{nmm} to this situation; see \refsec{whim}
for details.

The expression $E_f(p) = \psi_p \circ N_{f \circ \phi_p}(0_p)$
``lifts'' the Newton iteration function $N$ from Euclidean space
to a manifold.  \refsec{it} explores in generality the lifting of
an iteration function from Euclidean space to a manifold.

\subsection{Implications, Limitations and Examples}
\label{sec:imp}

Two broad types of optimisation problems can be distinguished.  One
is when little is known in advance about the possible cost functions
(save perhaps that they are convex, for example) and an algorithm
is desired that scales well with increasing dimension.  The other
is when the family of possible cost functions is known in advance
and an algorithm is desired that works well for all members of the
family.  The latter is the implicit focus of the current paper and
relates to real-time optimisation problems in signal processing:
at each instance, a new observation $y$ is made; this serves to
select a cost function $f(\cdot;y)$; it is required to find quickly
an $x$ that maximises $f(x;y)$.

Although generic choices are possible of the pair of parametrisations
$\phi$ and $\psi$ defining a Newton method on a particular manifold,
the fact remains that for large-scale problems, the Newton method
is generally abandoned in favour of quasi-Newton methods that build
up approximations of the Hessian over time, thereby making computational
savings by not evaluating the Hessian at each iteration.  Quasi-Newton
methods have memory and thus are not of the form $x_{k+1} = N_f(x_k)$.
An intended sequel will study how to lift algorithms with memory
to manifolds.

\emph{How can a Newton method be customised for a given family of
cost functions?}  It is propounded that thinking in terms of
parametrisations $\phi$ and $\psi$ offers greater insight into the
design of optimisation algorithms.  Notwithstanding that identifying
a ``killer application'' for the theory is work in progress, the
following example may sway some readers.

Generalising the Rayleigh quotient to higher dimensions yields two
well-studied optimisation problems~\cite{Helmke:1994ec}.  Recall
that the $(n,p)$-Stiefel manifold is the set of matrices $X \in
\reals^{n \times p}$ satisfying $X^T X = I$, where superscript $T$
denotes transpose and $I$ is the identity matrix.  (The manifold
structure is inherited from $\reals^{n \times p}$.) Let $A \in
\reals^{n \times n}$ be symmetric and $N \in \reals^{p \times p}$
diagonal, both with distinct positive eigenvalues.  A minimising
$X \in \reals^{n \times p}$ of $f(X) = \Tr(X^T A X N)$ subject to
$X^T X = I$ has as its columns the eigenvectors of $A$ corresponding
to the $p$ smallest eigenvalues of $A$.  If it was only required
to find the subspace spanned by these minor eigenvectors, known as
the minor subspace of $A$, then it suffices to minimise $g(X) =
\Tr(X^T A X)$ on the Grassmann manifold.  The Grassmann manifold
is a quotient space obtained from the Stiefel manifold by declaring
two matrices $X,Y$ as equivalent whenever there exists an orthogonal
matrix $Q$ such that $Y=XQ$.  In other words, each point on the
$(n,p)$-Grassmann manifold represents a particular $p$-dimensional
subspace of $\reals^n$.

The rate of convergence of $X_{k+1} = E_f(X_k) = \psi_{X_k} \circ
N_{f \circ \phi_{X_k}}(0_{X_k})$ is dictated by how close to being
quadratic $f \circ \phi_X$ is about $0_X$ whenever $X$ is near a
critical point of $f$.  As $f(X) = \Tr(X^T A X N)$ is already
quadratic, the parametrisation $\phi_X$ should be as linear as
possible.  One possibility is defining $\phi_X(Z)$ as the point on
the Stiefel manifold closest (in the Euclidean metric) to the matrix
$Z$; \refsec{proj} proves that parametrisations based on projections
are linear to at least second order.  The role of $\psi$ is to map
$N_{f \circ \phi_X}(0_X)$ back to the manifold with a minimum of
fuss.  Choosing $\psi$ to be the same as $\phi$ suffices.
(The option exists of choosing $\psi$ to be an approximation
of $\phi$ that makes $E_f(X) = \psi_X \circ N_{f \circ \phi_X}(0_X)$
overall less computationally demanding to evaluate numerically than
$E_f(X) = \phi_X \circ N_{f \circ \phi_X}(0_X)$.)  Since the Grassmann
manifold is a quotient of the Stiefel manifold, the above argument
readily extends to minimising $g(X)$ on a Grassmann manifold;
see~\cite{Manton:opt_mfold} for the precise calculations.

The above algorithm was trivial to derive yet is a sound starting
point upon which clever refinements are possible~\cite{Absil:2002js,
Absil:2004bd}.  The cubic rate of convergence is readily explained
in terms of $g \circ \phi_X$ being quadratic to third order at
critical points; compare with \refex{cubic}.  A feature of the
derivation is choosing $\phi_X$ with purpose rather than by trial
and error.

\emph{How should a theory of optimisation on manifolds be framed?}
This third italicised objective of the paper is in response to
misconceptions including: a connection is required for a Newton
method to be definable; only Riemannian Newton methods are ``true''
Newton methods; and, methods not exploiting the curvature of the
manifold must be inferior.  These misconceptions come from overplaying
the geometry of the manifold itself.

The most relevant geometry is that of the family of cost
functions~\cite{Manton:OG}.  Knowing the possible cost functions
$f$ allows for the customisation of the Newton method by choosing
a parametrisation $\phi$ that makes $f \circ \phi$ approximately
quadratic, and such a choice depends not on the manifold $M$ but
on the family of cost functions.  (Placing a sensible geometry on
$M$ might be advantageous --- perhaps computational burden can be
reduced by exploiting symmetry --- but the overall benefit nevertheless
will depend on the cost functions.)

It is not pragmatic to insist that only Riemannian Newton methods
\refeq{riem_n} are true Newton methods.  Different methods work
better for some cost functions and worse for others; no single
method can be superior for every smooth cost function. Any method
achieving a locally quadratic rate of convergence is worthy of the
title Newton method, provided of course it depends only on the 2-jet
of the function; see the definition given earlier.

Under this more general definition, there are Newton methods that
cannot be defined in terms of a connection.  A connection must vary
smoothly whereas no such requirement exists for the parametrisation
$\phi$.  More importantly, thinking of parametrisations instead of
connections is more conducive to customising a Newton method for a
given family of cost functions.  (The Riemannian approach \refeq{riem_n}
does not offer explicit insight into which metric to use if there
are two or more competing metrics, or what to do if there is no
convenient choice of metric.)

This paper avoids any need of Riemannian geometry by framing the
theory of optimisation on manifolds in terms of robustness of the
iteration function to changes of coordinates; see \refsec{it} for
details. This appears to be the most natural point of view.

\subsection{Motivation and Relationship with Other Work}

Given the extensive background and bibliography made available in
the book~\cite{Absil:2008uy}, only a handful of papers are discussed
below.

The Riemannian Newton method \refeq{riem_n} was introduced
in~\cite{Gabay:optm} but apparently went unnoticed.  The same
methodology was rediscovered in the influential
paper~\cite{Edelman:1998ei}.  The mindset is that the Newton method
is defined by its formula \refeq{newt}, and its extension to a
manifold thus necessitates endowing the manifold with a Riemannian
metric so the gradient and Hessian can be defined.

Numerically evaluating the Riemannian exponential map in \refeq{riem_n}
can be costly. It is common to replace the exponential map by an
approximation that is cheaper to evaluate numerically.  This is
formalised in~\cite{Adler:2002cc}, with a precursor in~\cite{Shub:1986ub}.
It corresponds to $E_f(p) = \psi_p \circ N_{f \circ \Exp_p}(0_p)$.
A Riemannian metric is still required for computing the Newton
increment, with what is termed a \emph{retraction} $\psi$ mapping the
result back to the manifold.  The retraction $\psi$ must satisfy
several conditions, including being smooth.  In the present paper,
$\psi$ need not be continuous and hence is not even a retraction
in the topological sense.  The way retractions are commonly used
in topology differs in spirit from how $\phi$ and $\psi$ are being
used to lift the Newton method to a manifold, hence the persistence
here of calling them parametrisations.

The Riemannian mindset was challenged in~\cite{Manton:opt_mfold}.
The basic idea is that since the Newton method is a local method,
the cost function in a neighbourhood of the current point can be
pulled back to a cost function on Euclidean space via a parametrisation,
one step of the Newton method carried out in Euclidean space, and
the result mapped back to the manifold.  No Riemannian metric is
necessary.  This corresponds to $E_f(p) = \phi_p \circ N_{f \circ
\phi_p}(0_p)$.  Using projections to define the parametrisations
$\phi_p$ was emphasised.  (The resulting algorithms differ significantly
from projected Newton methods that take a Newton step in the ambient
space then project back to the constraint surface.)

Combining the use of $\psi$ in~\cite{Adler:2002cc} and the use of
$\phi$ in~\cite{Manton:opt_mfold} immediately yields the general
form $E_f(p) = \psi_p \circ N_{f \circ \phi_p}(0_p)$ that is the
protagonist of the present paper.  This form is developed systematically
in Sections~\ref{sec:cad} and~\ref{sec:nmm}
in a way that suggests it is the most general form possible of a
Newton method.

The use of projections to define parametrisations, advocated
in~\cite{Manton:opt_mfold}, was studied in~\cite{Absil:2012eb},
but for $E_f(p) = \psi_p \circ N_{f \circ \Exp_p}(0_p)$.  Convergence
proofs were based on calculus techniques requiring more orders
of differentiability than necessary; see \refsec{NewtEucl}.

Another active stream of research is finding lower bounds on the
radius of convergence of Riemannian Newton methods~\cite{Alvarez:2008hv,
Argyros:2007cw, Ferreira:2002uj}.  This has not been addressed in
the present paper, although in principle, a careful study of the
constants in the bounds derived here would provide that information.

The question of the most general form of a Newton method on a
manifold appears not to have been addressed before.


\section{Basic Notation and Definitions}
\label{sec:rate_c}

For a function $f$ between Euclidean spaces, the following definitions
are made.  The Euclidean norm $\| \cdot \|$ on $\rn$ is used
throughout.  The norm of the second-order derivative $D^2f(x)$ is
$\| D^2f(x) \| = \sup_{\|\eta\|=1} \|D^2f(x)\cdot(\eta,\eta)\|$.
All other norms are operator norms.  Gradients $\nabla f$ and
Hessians $\Hess_f$ are calculated with respect to the Euclidean
inner product.  The identity operator is denoted by $I$ (or sometimes
by $1$ in the one-dimensional case).  The notation $B_n(x;\rho)$ and its
abbreviation $B(x;\rho)$ denote the
open ball centred at $x \in \rn$ of radius $\rho$. Its closure is
$\overline{B}(x;\rho)$.

An iteration function $N\colon \rn \rightarrow \rn$, which may not be defined
on the whole of $\rn$, is said to converge locally to $\xstar \in \rn$
with rate $K \in \reals$ and constant $\kappa \in \reals$ if there
exists an open set $U \subset \rn$ containing $\xstar$ such that $N$
is defined on $U$ and
\begin{equation}
\label{eq:cnv}
\forall x \in U,\quad N(x) \in U\ \mathrm{and}\ 
\| N(x) - \xstar \| \leq \kappa\,\| x - \xstar \|^K.
\end{equation}
If $K=1$ then it is further required that $\kappa < 1$,
and convergence is called linear.  If $K \in (1,2)$ the convergence
is super-linear, and if $K=2$ the convergence is quadratic.

Although \refeq{cnv} implies $N(\xstar) = \xstar$, the sequence $x_{k+1}
= N(x_k)$ need not converge to $\xstar$ for an arbitrary $x_0 \in U$.
Nevertheless, define $\bar{\rho} = \kappa^{1/(1-K)}$ if $K>1$, or
$\bar{\rho} = \infty$ if $K=1$ and $\kappa < 1$.  Then $B(\xstar;\rho)
\cap U$ is mapped into itself by $N$ whenever $\rho \leq \bar{\rho}$.
Moreover, $x_0 \in B(\xstar;\bar{\rho}) \cap U$ implies $x_k \rightarrow
\xstar$.

The focus of this paper is on convergence rates $K$ greater than one.


\section{Local Convergence of the Newton Iteration on Euclidean Space}
\label{sec:NewtEucl}

Convergence proofs for the Newton method include the Newton-Kantorovich
theorem (applicable for the Newton method on Banach spaces) and the
Newton-Mysovskikh theorem;
see~\cite{bk:Kantorovich:fn_analysis,bk:Ortega:iter_soln} and the
bibliographic note~\cite[p.\@ 428]{bk:Ortega:iter_soln}.  These
theorems give sufficient but not necessary conditions, concentrating
instead on explicitly finding a region within which the Newton
method is guaranteed to converge.  The affine invariance of the
Newton method is exploited in~\cite{bk:Deuflhard:newton} to sharpen
these classical results.

In pursuit of the most general Newton method on a manifold, it is
informative to derive a necessary and sufficient condition for the
standard Newton method to converge to a non-degenerate critical
point.  

\begin{thm}
\label{th:Neuclid}
Let $f\colon \rn \rightarrow \reals$ be $\Ct^2$-smooth.  Let $\xstar \in
\rn$ be a non-degenerate critical point, that is, $\nabla f(\xstar)=0$ and
$\Hess_f(\xstar)$ is invertible.  A necessary and sufficient condition for
$N_f$ in \refeq{newt} to be locally quadratically convergent to $\xstar$
is for there to exist $\eta,\delta > 0$ such that $x \in B(\xstar;\delta)$
implies
\begin{equation}
\label{eq:Neuclid}
\| \left[\Hess_f(x)-\Hess_f(\xstar)\right] (x-\xstar) \|
	\leq \eta\,\|x-\xstar\|^2.
\end{equation}
\end{thm}
\begin{proof}
Define the second-order Taylor series remainder term
\begin{equation}
\label{eq:Rx}
R(x) = f(x) - f(\xstar) - \frac12 (x-\xstar)^T \Hess_f(\xstar) (x-\xstar).
\end{equation}
Since $f$ is $\Ct^2$, so is $R$.  Moreover,
\begin{align}
\nabla f(x) &= \Hess_f(\xstar) (x-\xstar) + \nabla R(x), \\
\Hess_f(x) &= \Hess_f(\xstar) + \Hess_R(x).
\end{align}
Substitution into \refeq{newt} shows
\begin{align}
N_f(x)-\xstar &= x - \xstar - \left[\Hess_f(x) \right]^{-1}
	\left[\Hess_f(\xstar) (x-\xstar) + \nabla R(x)\right] \\
\label{eq:NffR}
  &= \left[\Hess_f(x) \right]^{-1}
    \left[ \Hess_R(x) (x-\xstar) - \nabla R(x) \right].
\end{align}
Since $\Hess_R$ is continuous, for any $\epsilon > 0$
there exists a $\rho > 0$ such that $x \in B(\xstar;\rho)$ implies:
$\Hess_f(x)$ is invertible; $\left\| \left[ \Hess_f(x) \right]^{-1}
\right\| \leq \left\| \left[ \Hess_f(\xstar) \right]^{-1} \right\|
+ \epsilon$; and $\left\| \Hess_f(x) \right\| \leq \left\|
\Hess_f(\xstar) \right\| + \epsilon$.  

To prove sufficiency, first observe
\begin{equation}
\| N_f(x) - \xstar \| \leq 
	\left\| \left[\Hess_f(x) \right]^{-1} \right\|
	\left(\| \Hess_R(x) (x-\xstar)\| + \| \nabla R(x) \|\right).
\end{equation}
Choose $\delta,\eta$ as in the theorem.  If
$x \in B(\xstar;\delta)$ then
\begin{align}
\| \nabla R(x) \|
  &\leq \int_0^1 \frac1t \| \Hess_R(\xstar + t(x-\xstar)) t(x-\xstar) \|\,dt \\
  &\leq \int_0^1 t\,\eta \| x - \xstar \|^2\,dt \\
\label{eq:nabrx}
	&\leq \frac12 \eta\,\| x - \xstar \|^2.
\end{align}
Choosing $\epsilon,\rho$ as above,
if $x \in B(\xstar;\min\{\delta,\rho\})$ then $N_f(x)$ is well-defined and
\begin{equation}
\| N_f(x) - \xstar \| \leq \frac32 \eta
	\left( \left\| [\Hess_f(\xstar)]^{-1} \right\| + \epsilon \right)
	\| x - \xstar \|^2,
\end{equation}
proving local quadratic convergence.

To prove necessity, first note from \refeq{NffR} that
\begin{equation}
\| N_f(x) - \xstar \| \geq
	\left\| \Hess_f(x) \right\|^{-1}
	    \left\| \Hess_R(x) (x-\xstar) - \nabla R(x) \right\|.
\end{equation}
Thus, choosing $\epsilon,\rho$ as above,
if $x \in B(\xstar;\rho)$ then
\begin{equation}
\left\| \Hess_R(x) (x-\xstar) - \nabla R(x) \right\| \leq
\left( \left\| \Hess_f(\xstar) \right\| + \epsilon \right)
\| N_f(x) - \xstar \|.
\end{equation}
By hypothesis, $N_f$ converges locally quadratically to $\xstar$, hence
by shrinking $\rho$ if necessary, there exists a $\kappa > 0$ such
that $x \in B(\xstar;\rho)$ implies
\begin{equation}
\label{eq:Rkappa}
\left\| \Hess_R(x) (x-\xstar) - \nabla R(x) \right\|
	\leq \kappa \| x - \xstar \|^2.
\end{equation}
Define the closed ball $C = \overline{B}(\xstar;\rho/2)$ and the
function $\phi(x) = \| \Hess_R(x) (x-\xstar) \| \|x-\xstar\|^{-1}$.
Setting $\phi(\xstar)=0$ ensures $\phi$ is well-defined and continuous
on $C$.  Assume to the contrary, for all $\eta > 0$, the scalar
$h = \max_{x \in C} \left\{ \phi(x) - \eta\|x-\xstar\| \right\}$
satisfies $h>0$.  For any $x \in C$,
\begin{align}
\| \nabla R(x) \|
  &\leq \int_0^1 \| \Hess_R(\xstar+t(x-\xstar)) (x-\xstar)\|\,dt \\
  &= \|x-\xstar\| \int_0^1 \phi(\xstar+t(x-\xstar))\,dt \\
  &\leq \|x-\xstar\| \int_0^1 h + t\eta\|x-\xstar\|\,dt \\
  &= h \|x-\xstar\| + \frac12\eta\|x-\xstar\|^2.
\end{align}
Let $z \in C$ be such that $\phi(z)-\eta\|z-\xstar\|=h$.
Since $z \neq \xstar$ and
\begin{align}
\| \Hess_R(z)(z-\xstar) - \nabla R(z)\| &\geq \phi(z) \|z-\xstar\|
  - h \|z-\xstar\| - \frac12\eta\|z-\xstar\|^2 \\
\label{eq:Reta}
  &= \frac12\eta\|z-\xstar\|^2,
\end{align}
choosing $\eta > 2\kappa$ makes \refeq{Reta} contradict \refeq{Rkappa},
proving the theorem.
\end{proof}

\begin{cor}
Let $f\colon \rn \rightarrow \reals$ be $\Ct^3$-smooth and $\xstar \in \rn$
a non-degenerate critical point.  Then $N_f$ in \refeq{newt} converges
locally quadratically to $\xstar$.
\end{cor}
\begin{proof}
If $f$ is $\Ct^3$ then $\Hess_f(x)-\Hess_f(\xstar)$ is $\Ct^1$,
hence \refeq{Neuclid} holds.
\end{proof}

\begin{cor}
\label{cor:modN}
Let $f$ and $\xstar$ satisfy the conditions in \refth{Neuclid},
including \refeq{Neuclid}.  The perturbed iteration function $E_f(x) = x -
[\Hess_f(x) + G(x)]^{-1} \nabla f(x)$ converges locally quadratically
to $\xstar$ if there exists a $\gamma \in \reals$ such that the
operator norm of the
matrix $G(x)$ satisfies $\| G(x) \| \leq \gamma \| x - \xstar \|$
in a neighbourhood of $\xstar$.
\end{cor}
\begin{proof}
Observe
\begin{equation}
E_f(x)-\xstar = [\Hess_f(x)+G(x)]^{-1}\left\{
\Hess_f(x)(x-\xstar) - \nabla f(x) + G(x)(x-\xstar) \right\}.
\end{equation}
Therefore,
\begin{multline}
\| E_f(x) - \xstar \|
  \leq \left\| [\Hess_f(x)+G(x)]^{-1} \right\| \big\{
  	\| \nabla f(x) - \Hess_f(\xstar)(x-\xstar) \| + \\
	\left\| [\Hess_f(x) - \Hess_f(\xstar)](x-\xstar) \right\| +
	\| G(x) (x-\xstar) \| \big\}.
\end{multline}
In a sufficiently small neighbourhood of $\xstar$, $\left\|
[\Hess_f(x)+G(x)]^{-1} \right\|$ is bounded above by a constant and
the three other terms are bounded by a constant times $\|x-\xstar\|^2$;
refer to \refeq{nabrx} and the hypotheses on $\Hess_f(x)$ and $G(x)$.
\end{proof}

Despite calculus offering a simpler and more elegant alternate,
convergence proofs are based here on hard analysis because calculus
requires a higher order of smoothness than necessary, as now
demonstrated.  (See also the opening paragraph of \refapp{sc}.)
Recall the basic principle.

\begin{lem}
\label{lem:Sdiff}
Let $N\colon \rn \rightarrow \rn$ be $\Ct^K$-smooth for some integer $K
\geq 2$.  If $D^k N(\xstar) = 0$ for $k=1,\cdots,K-1$ then $N$
converges locally to $\xstar$ with rate $K$.
\end{lem}

Applying \reflem{Sdiff} to \refeq{newt} shows that $f$ being
$\Ct^4$-smooth is sufficient for $N_f$ to converge locally quadratically
to a nondegenerate critical point.  If $f$ were only $\Ct^3$ then
$N_f$ would only be $\Ct^1$ and \reflem{Sdiff} could not be applied.
The actual condition \refeq{Neuclid} falls strictly between
$\Ct^2$-smoothness and $\Ct^3$-smoothness.

\begin{exmp}
Define $f(x) = x^2 + | x |^{5/2}$.  The origin is a non-degenerate
critical point.  The Newton iteration function is $N_f(x) = \frac{5x
|x|^{1/2}}{8+15 |x|^{1/2}}$ and has super-linear but not quadratic
convergence, despite $f$ being $\Ct^2$-smooth.
\end{exmp}

\begin{rem}
\label{rem:Neuclidind}
The quadratic convergence rate of the Newton method
is coordinate independent, in the following sense.  Assume $f\colon \rn
\rightarrow \reals$ satisfies the conditions in \refth{Neuclid}
about the point $\xstar$. If $\phi\colon \rn \rightarrow \rn$ is a
$\Ct^2$-diffeomorphism then $\phi^{-1}(\xstar)$ is a non-degenerate
critical point of $f \circ \phi$, and by \refprop{intHess}, condition
\refeq{Neuclid} holds for $f \circ \phi$ about the point
$\phi^{-1}(\xstar)$.  Thus, if $N_f$ converges locally quadratically
to $\xstar$ then $N_{f \circ \phi}$ converges locally quadratically
to $\phi^{-1}(\xstar)$.
\end{rem}


\section{The Coordinate Adapted Newton Iteration}
\label{sec:cad}

The most general form of a Newton method in Euclidean space is
explored.

\subsection{Coordinate Adaptation}
\label{sec:canewt}

Applying a change of coordinates $\phi\colon \rn \rightarrow \rn$ to
\refeq{newt} yields the new iteration function $E_f(x) = \phi
\circ N_{f \circ \phi} \circ \phi^{-1}(x)$.  Expedient choices of
$\phi$ can increase the domain of attraction, decrease the computational
complexity per iteration and improve the convergence rate.  As an
extreme example, if $\phi$ is such that $f \circ \phi$ is quadratic
then $E_f$ converges in a single iteration.  Although Morse's Lemma
guarantees the existence of such a $\phi$ locally, finding it is
generally not practical.  This motivates using a different change
of coordinates at each iteration, namely $E_f(x) = \phi_x \circ
N_{f \circ \phi_x} \circ \phi_x^{-1}(x)$.  When $\phi_x$ varies
with $x$, the convergence properties of $E_f$ need not 
follow from the convergence properties of $N_f$. Significantly then,
it is established that under mild conditions, $E_f$ converges locally
quadratically to non-degenerate critical points of $f$.

Coordinate adaptation is defined in terms of a function
$\phi\colon \rn \times \rn \rightarrow \rn$, alternatively written
$\phi_x(y) = \phi(x,y)$, satisfying the condition
that, $\forall \xstar \in \rn$, $\exists \alpha, \beta, \rho \in
\reals$, $\rho > 0$, $\forall x,y \in B(\xstar;\rho)$,
the following hold:
\begin{description}
\item[P1] $\| D^2\phi_x(x) \| \leq \alpha$;
\item[P2] $\| \phi_x(y) - y \| \leq \beta \| y-x \|^2$.
\end{description}
Implicit in P1 is the requirement that $D^2\phi_x(x)$ exists, which
in turn requires the existence of $D\phi_x(y)$ for $y$ sufficiently
close to $x$.  

Given such a $\phi$, the coordinate adapted Newton iteration function is
\begin{equation}
\label{eq:Ef}
E_f(x) = \phi_x \circ N_{f \circ \phi_x}(x)
\end{equation}
where $N_f$ is the Newton iteration function \refeq{newt}.  This agrees
with the earlier expression for $E_f$ because P2 implies
$\phi_x(x)=x$.

\begin{thm} \label{th:cad}
Let $f$ and $\xstar$ satisfy the conditions in \refth{Neuclid},
including \refeq{Neuclid}.  Let $\phi$ satisfy P1 and P2, defined
above.  Then the coordinate adapted Newton iteration function $E_f$, defined
in \refeq{Ef}, converges locally quadratically to $\xstar$.
\end{thm}
\begin{proof}
P2 implies $\phi_x(x)=x$ and $D\phi_x(x)=I$.  Hence
\begin{align}
D(f \circ \phi_x)(x) &= Df(x), \\
D^2(f \circ \phi_x)(x) &= D^2f(x) + Df(x) D^2\phi_x(x).
\end{align}
Let $G(x)$ be the matrix representation of $Df(x) D^2\phi_x(x)$.
Then $G(x)$ is symmetric and satisfies $\xi^T G(x) \xi = Df(x)
D^2\phi_x(x) \cdot (\xi,\xi)$ for any $\xi \in \rn$.  If P1 holds,
it can be shown that $\| G(x) \| \leq \alpha \| Df(x) \|$.  Thus,
in a neighbourhood of $\xstar$, there exist constants $r,\gamma >
0$ such that $x \in B(\xstar;r)$ implies $\| G(x) \| \leq \gamma
\|x-\xstar\|$.  Since $\Hess_{f \circ \phi_x}(x) = \Hess_f(x) +
G(x)$, it follows from \refcor{modN} that, for a possibly smaller
$r>0$, there exists a $\kappa > 0$ such that $\| N_{f \circ \phi_x}(x)
- \xstar \| \leq \kappa \| x-\xstar\|^2$ whenever $x \in B(\xstar;r)$.
To be able to apply P2, shrink $r$ if necessary to ensure $0 < r
\leq \rho$.  Then
\begin{align}
\| E_f(x) - \xstar \|
  &\leq \| \phi_x\left( N_{f \circ \phi_x}(x) \right) -
  	N_{f \circ \phi_x}(x) \| + \| N_{f \circ \phi_x}(x) - \xstar \| \\
  &\leq \beta \| N_{f \circ \phi_x}(x) - x \|^2 + \kappa \| x-\xstar\|^2 \\
  &\leq \beta \left( \| N_{f \circ \phi_x}(x) - \xstar \| +
  	\| x-\xstar\| \right)^2 + \kappa \| x-\xstar\|^2 \\
  &\leq \left( \beta (\kappa r + 1)^2 + \kappa \right) \|x-\xstar\|^2
\end{align}
whenever $x \in B(\xstar;r)$, proving the theorem.
\end{proof}

As now explained, P1 and P2 are not only mild, it is conjectured
they cannot be weakened.  For \refeq{Ef} to be defined, the Hessian
of $f \circ \phi_x$ at $x$ must exist, necessitating the existence
of $D^2\phi_x(x)$ implicit in P1.  The local bound in P1 ensures
$N_{f \circ \phi_x}(x)$ converges locally quadratically.  A side-effect
of P2 is that $\phi_x(x) = x$ and $D\phi_x(x)=I$;
this loses no generality because the Newton iteration
function is invariant to affine changes of coordinates.  The main
purpose of P2 is to prevent the residual term $R_x(y) = \phi_x(y)
- y - \frac12 D^2\phi_x(x) \cdot (y-x,y-x)$ from being unbounded
locally, ensuring that if $E(x)$ is an arbitrary iteration function
converging locally quadratically to $\xstar$ then $\phi_x \circ
E(x)$ continues to converge locally quadratically to $\xstar$.  The
situation in which, for a sufficiently large class of cost functions
$f$, $N_{f \circ \phi_x}(x)$ fails to have local quadratic convergence
yet $E_f$ has local quadratic convergence is conjectured to be
impossible.  The claim that P1 and P2 are mild comes from the fact
deducible from \refapp{cond}
that any $\Ct^2$-smooth $\phi$ with $\phi_x(x) =
x$ and $D\phi_x(x)=I$, satisfies P1 and
P2.  Furthermore, neither $\phi(x,y)$ or $D\phi_x(y)$ need be
continuous except on the diagonal $y=x$, and $\phi_x$ need not be
locally $\Ct^1$-smooth.

The following example shows that if arbitrary changes of coordinates
are allowed then the coordinate adapted Newton method may not even
be defined, much less converge at a quadratic rate.

\begin{exmp}\label{ex:beta}
Let $\beta$ be an arbitrary scalar.  Consider the coordinate adapted
Newton iteration function applied to $f(x) = x^2$ using $\phi_x(y) = y +
\frac\beta x(y-x)^2$ when $x \neq 0$ and $\phi_x(y) = y$ when $x=0$.
If $x \neq 0$ then $N_{f \circ \phi_x}(x) = \frac{2\beta}{1+2\beta}x$
which is not defined if $\beta = -\frac12$.  If $\beta \neq -\frac12$
then $E_f(x) = \frac{(3+4\beta)\beta}{(1+2\beta)^2} x$, which in
general exhibits at best linear convergence.
\end{exmp}

\subsection{The Generalised Coordinate Adapted Newton Iteration}
\label{sec:gcanewt}

The proof of \refth{cad} can be modified trivially to prove the
following result.  A new function $\psi$ analogous to $\phi$ is
introduced and property P2 in \refsec{canewt} is replaced by
\begin{description}
\item[P2'] $\| \psi_x(y) - y \| \leq \beta \| y-x \|^2$.
\end{description}

\begin{thm} \label{th:gcad}
Let $f$ and $\xstar$ satisfy the conditions in \refth{Neuclid},
including \refeq{Neuclid}.  Let $\phi\colon \rn \times \rn \rightarrow
\rn$ satisfy P1 in \refsec{canewt}.  Assume further that
$\phi_x(x) = x$ and $D\phi_x(x)=I$.  Let $\psi\colon \rn
\times \rn \rightarrow \rn$ satisfy P2' above; the qualifiers for
$x$ and $y$ in P2' are the same as for P2.
Then the generalised coordinate adapted Newton iteration function
\begin{equation} \label{eq:gcad}
E_f(x) = \psi_x \circ N_{f \circ \phi_x}(x)
\end{equation}
converges locally quadratically to $\xstar$.
\end{thm}

Being able to change both $\phi$ and $\psi$ allows greater control
over the computational complexity, the domain of attraction and the
rate of convergence of the iteration function \refeq{gcad}, as now discussed.

\subsection{Discussion}
\label{sec:cad_disc}

The choice of coordinate changes $\phi_x$ and $\psi_x$ in \refeq{gcad}
determines which class of cost functions the generalised coordinate
adapted Newton method will perform well for.  The challenge then
is to determine suitable coordinate changes to use for the class
of cost functions at hand.  For inherently difficult optimisation
problems this will not be easy by definition.  Nevertheless, thinking
in terms of coordinate adaptation leads to the following new strategy.

The closer the cost function is to being quadratic, the faster the
convergence rate of \refeq{newt}.  Ideally then, $\phi_x$ in
\refeq{Ef} makes $f \circ \phi_x$ approximately quadratic for every
cost function $f$ in the given family.  For improving local
convergence, it suffices to restrict attention to cost functions
with a critical point near $x$ because, by definition of local, it
can be assumed a critical point is nearby $x$, and this limits the
possibilities of which cost function from the family has been
selected to be minimised.  For the special case when $f$ has the
form $f(x) = \psi(x-z)$ for some unknown scalar $z$, where $\psi$
has a minimum at the origin, it suffices for $f \circ \phi_x$ to
be approximately quadratic when $z=x$.

\begin{exmp}
\label{ex:cubic}
Consider the family of cost functions $f(x;z) = (x-z)^2 + 2(x-z)^3$.
The coordinate adapted Newton iteration function using the coordinate
systems $\phi_x(y) = y - (y-x)^2$ is $E_f(x;z) = z - 8(x-z)^3 +
\cdots$.  This converges cubically to the critical
point $\xstar = z$ for any cost function in the family.  Here,
$\phi_x(y)$ was chosen so $f \circ \phi_z(x) = (x-z)^2
- 5 (x-z)^4 + \cdots$ has no cubic term.
\end{exmp}

If the domain of attraction is of primary concern then similar
intuition suggests choosing $\phi_x$ such that, for any $f$ belonging
to the given class of cost functions, $f \circ \phi_x$ has a
relatively large domain of attraction, especially if $x$ is at all
close to the minimum of $f$.

The extra freedom afforded by $\psi_x$ in \refeq{gcad} can be used
to reduce the computational complexity per iteration without
compromising the rate of convergence; in some cases, an expedient
choice of $\psi_x$ leads to cancellations, so $\psi_x \circ N_{f
\circ \phi_x}$ becomes less computationally intensive to evaluate
than $N_{f \circ \phi_x}$ on its own.

The coordinate adapted Newton method is different from variable
metric methods.  Variable metric methods explicitly or implicitly
perform a change of coordinates and then take a steepest-descent
(not Newton)
step in the new coordinate system.  They do not evaluate the Hessian
of the cost function but instead build up an approximation $B_k$
to the Hessian from current and past gradient information.  They
are of the form $x_{k+1} = x_k - [B_k]^{-1} \nabla f(x_k)$;
see~\cite{bk:Polak:opt}.
The generalised coordinate adapted Newton method \refeq{gcad} with
$\psi_x(y)=y$ can be written as $E_f(x) = x - [\Hess_f(x)+G(x)]^{-1}
\nabla f(x)$; see the proof of \refth{cad}.  This differs from a
variable metric method in several ways; $E_f$ makes use of the
Hessian of $f$ but $B_k$ does not; $E_f$ has no ``memory'' but $B_k$
is built up over time; variable metric methods generally only achieve
super-linear convergence whereas $E_f$ has quadratic convergence.
The philosophy is also different; variable metric methods wish for
$B_k$ to be as close as possible to the true Hessian, whereas the
generalised coordinate adapted Newton method intentionally uses a
perturbed version of the true Hessian to improve the performance
of the algorithm.


\section{Generalised Newton Methods on Manifolds}
\label{sec:nmm}

Throughout this section, $f\colon M \rightarrow \reals$ will be a
$\Ct^2$-smooth cost function defined on an $n$-dimensional
$\Ct^2$-differentiable manifold $M$.  Recall from \refapp{conv}
that if $E\colon M \rightarrow M$ is an iteration function with local
quadratic convergence to a point $\pstar \in M$ then $\varphi \circ
E \circ \varphi^{-1}$ converges locally quadratically to $\varphi(\pstar)$
for any chart $(U,\varphi)$ with $\pstar \in U$.  Fix $(U,\varphi)$.
For cost functions $f$ with a critical point in $U$, the coordinate
adapted Newton method of \refsec{cad} can be extended to manifolds
by seeking an $E_f$ such that $\varphi \circ E_f \circ \varphi^{-1}$
is a coordinate adapted Newton iteration function for the equivalent
cost function $f \circ \varphi^{-1}$.  For functions with critical
points outside $U$, in principle a different coordinate chart
needs to be taken, but as shown presently, it is straightforward
to guess an appropriate form for $E_f$ globally.

Solving $\varphi \circ E_f \circ \varphi^{-1}(x)
= \phi_x \circ N_{\left(f \circ \varphi^{-1}\right) \circ \phi_x}
\circ \phi_x^{-1}(x)$ yields
\begin{align}
E_f(p) &= \left(\varphi^{-1} \circ \phi_{\varphi(p)}\right) \circ
    N_{f \circ \left(\varphi^{-1} \circ \phi_{\varphi(p)}\right)} \circ
    \left(\phi_{\varphi(p)}^{-1} \circ \varphi\right)(p) \\
    &= \psi_p \circ N_{f \circ \psi_p} \circ \psi_p^{-1}(p), \qquad
        \psi_p(z) = \varphi^{-1} \circ \phi_{\varphi(p)}(z).
\end{align}
The affine invariance of the Newton method allows this to be rewritten
as $E_f(p) = \tilde\psi_p \circ N_{f \circ \tilde\psi_p} \circ
\tilde\psi_p^{-1}(p) = \tilde\psi_p \circ N_{f \circ \tilde\psi_p}(0)$
where $\tilde\psi_p(z) = \psi_p(z+\varphi(p))$.  Although this
defines $E_f$ only locally, an obvious extension is $E_f(p) = \psi_p
\circ N_{f \circ \psi_p}(0)$ where, for each $p \in M$, $\psi_p\colon
\reals^n \rightarrow M$ is a parametrisation of a neighbourhood on
the $n$-dimensional manifold $M$ centred at $p$, that is, $\psi_p(0)=p$.
This extension is justified by the proof of \refth{nmm} in which it
is shown that $\varphi \circ E_f \circ \varphi^{-1}$ does indeed
take the form of a coordinate adapted Newton method for any chart
$\varphi$ on $M$.

Although tempting to generalise $\phi$ and $\psi$ in \refsec{cad}
to maps from $M \times \reals^n$ to $M$, the global geometry of $M$
can prevent any such map from being smooth.  The tangent bundle
$TM$, being equivalent to $M \times \reals^n$ locally, offers an
alternative.  As $TM$ twists in the ``right'' way, smooth
parametrisations from $TM$ to $M$ can be anticipated to exist;
this was appreciated by Shub~\cite{Shub:1986ub, Adler:2002cc}.
(While smoothness is not essential, in practice it may be convenient
to work with smooth parametrisations.)

\newcommand{\lc}[1]{\frac\partial{\partial\varphi_{#1}}}
\newcommand{\lca}[2]{\left.\frac\partial{\partial\varphi_{#1}}\right|_{#2}}

The functions $\phi\colon TM \rightarrow M$ and $\psi\colon TM \rightarrow
M$ will be required to satisfy conditions C1--C2 below, which generalise
P1 and P2' in Sections~\ref{sec:canewt} and~\ref{sec:gcanewt}.
Local coordinates are needed.  Let $\pi\colon TM \rightarrow M$ be the
projection taking a tangent vector $v_p \in T_pM$ to its base
point $p$.  A $\Ct^2$-chart $(U,\varphi)$ induces the $\Ct^1$-chart
$\tau_\varphi\colon \pi^{-1}(U) \rightarrow \reals^n \times \reals^n$
on $TM$, sending $v_p$ to $(\varphi(p),A_p(v_p))$ where $A_p\colon T_pM
\rightarrow \reals^n$ is the linear isomorphism taking $v_p =
\sum_{j=1}^n \eta^j \lca{j}{p}$ to $\eta =
(\eta^1,\cdots,\eta^n)$.  The local coordinate representation of
$\phi$ is $\widehat\phi = \varphi \circ \phi \circ
\tau_\varphi^{-1}$.

Conditions C1--C2 are satisfied if,
$\forall p \in M$,
$\exists$ $\Ct^2$-chart $(U,\varphi)$ with $\varphi(p)=0$,
$\exists \rho > 0$:
\begin{description}
\item[C1] $\widehat\phi=\varphi \circ \phi \circ \tau_\varphi^{-1}$
satisfies \Hi{\rho} and \Hii{\rho}, defined below;
\item[C2] $\widehat\psi=\varphi \circ \psi \circ \tau_\varphi^{-1}$
satisfies \Hiii{\rho}, defined below.
\end{description}

Consider a function $h\colon \reals^n \times \reals^n \rightarrow
\reals^n$ which need not be defined everywhere; its domain of
definition will be clarified shortly. It satisfies \textbf{\Hi{\rho}}
if $x \in B(0;\rho)$ implies $h_x(0)=x$ and $Dh_x(0)=I$.  It satisfies
\textbf{\Hii{\rho}} if there exists a constant $\alpha \in \reals$
such that $\| D^2h_x(0) \| \leq \alpha$ for every $x \in B(0;\rho)$.
It satisfies \textbf{\Hiii{\rho}} if there exists a constant $\beta
\in \reals$ such that $x,y \in B(0;\rho)$ implies $\| h_x(y) - x -
y \| \leq \beta \|y\|^2$.  If the subscript $\rho$ is omitted, the
existence of an appropriate $\rho>0$ is implied.

For the derivatives to exist, if $h$ satisfies \Hi{\rho} or \Hii{\rho}
then its domain of definition must include a set of the form $\{(x,y)
\mid x \in B(0;\rho),\ y \in B(0;\delta_x),\ \delta_x > 0\}$ where
$\delta_x$ is a function of $x$. Such a set need not contain a
neighbourhood of the origin.  For \Hiii{\rho} though, it is required
that $B(0;\rho) \times B(0;\rho)$ lies in the domain of $h$.  See
\refapp{cond} for further properties.

A generalised Newton iteration function is any 
$E_f\colon M \rightarrow M$ of the form
\begin{equation}
\label{eq:gni}
E_f(p) = \psi_p \circ N_{f \circ \phi_p}(0_p)
\end{equation}
where $\phi_p\colon T_pM \rightarrow M$ and $\psi_p\colon T_pM \rightarrow
M$ are the restrictions of $\phi$ and $\psi$ to the tangent space
$T_pM$ at the point $p$ on $M$.  In \refeq{gni}, $N$ represents the
Newton iteration \refeq{newt} but on the abstract vector space $T_pM$
rather than $\reals^n$.

The local coordinate representation $\widehat\phi = \varphi \circ
\phi \circ \tau_\varphi^{-1}$ of $\phi$ can be written as
$\widehat\phi(x,y) = \varphi \circ \phi_{\varphi^{-1}(x)} \circ
A_{\varphi^{-1}(x)}^{-1}(y)$, alternatively denoted $\widehat\phi_x(y)$.
The local coordinate representation of $\phi_p$ is $\varphi \circ
\phi_p \circ A_p^{-1}$, which can be written in terms of $\widehat\phi$,
namely, $\widehat\phi_{\varphi(p)}$.  Analogously for $\psi_p$.

\begin{thm}
\label{th:nmm}
Let $f\colon M \rightarrow \reals$ be a $\Ct^2$-smooth cost function on
a $\Ct^2$-smooth manifold $M$.  Let $\pstar \in M$ be a non-degenerate
critical point, that is, $Df(\pstar)=0$ and if
$D^2f(\pstar)\cdot(\xi,\xi)=0$ then $\xi=0$.  Assume there is a
chart $(U,\varphi)$ on $M$, with $\varphi(\pstar)=0$, and an $\eta
> 0$ such that $\widehat f = f \circ \varphi^{-1}$ satisfies,
for $x \in \varphi(U)$,
\begin{equation}
\label{eq:nmmf}
\| [\Hess_{\widehat f}(x) - \Hess_{\widehat f}(0)] x \| \leq \eta\|x\|^2.
\end{equation}
Let $\phi,\psi\colon TM \rightarrow M$ satisfy C1--C2, defined above.
Then the generalised Newton iteration function \refeq{gni} converges
locally quadratically to $\pstar$.
\end{thm}
\begin{proof}
Let $(U,\varphi)$ be as in the theorem.  \refprop{Cinv} implies there
exists a $\rho>0$ such that $\widehat\phi=\varphi \circ \phi \circ
\tau_\varphi^{-1}$ satisfies \Hi{\rho} and \Hii{\rho}, and
$\widehat\psi=\varphi \circ \psi \circ \tau_\varphi^{-1}$ satisfies
\Hiii{\rho}.  (By \reflem{PP}, $\widehat\psi$ will also satisfy \Hi{\rho}.)
Let $h_x(y) = y-x$, $\tilde\phi_x(y)
= \widehat\phi(x,y-x)
= \widehat\phi_x \circ h_x(y)$ and $\tilde\psi_x(y)
= \widehat\psi(x,y-x)
= \widehat\psi_x \circ h_x(y)$.
The invariance of the Newton iteration function
to the affine coordinate change $h_x^{-1} \circ A_{\varphi^{-1}(x)}$
can be used to show 
\begin{equation}
\varphi \circ E_f \circ \varphi^{-1}(x) =
    \tilde\psi_x \circ N_{\widehat f \circ \tilde\phi_x}(x).
\end{equation}
The functions $\tilde\phi$, $\tilde\psi$ and $\widehat f$ satisfy
the necessary conditions locally about the point $\xstar=0$ for the
proof of \refth{gcad} to go through.
\end{proof}

If $f$ is $\Ct^3$-smooth then \refeq{nmmf} holds.
See also \refapp{fr}.

\subsection{A Global Topological Constraint}
\label{sec:gtc}

The generalised Newton method \refeq{gni} is defined by the choice
of parametrisations $\phi$ and $\psi$, and \reflem{Mthree} shows
it suffices to choose $\phi$ and $\psi$ to be $\Ct^2$-smooth to
ensure local quadratic convergence to non-degenerate critical points.
This may appear elegant and straightforward.  In practice though,
especially if $M$ is a quotient space, directly writing down a
smooth parametrisation $\phi$ may not be the most desirable approach.
While numerous possible choices may come to mind for each $\phi_p\colon
T_pM \rightarrow M$, difficulties arise if no \textit{canonical}
choice is evident for each $p$ that would make $\phi$ smooth.

This difficulty is a consequence of a deeper fact: a global topological
constraint is unwittingly imposed by insisting that a Newton method
be strictly of the form $p_{k+1} = E_f(p_k)$.  To see this, consider
assigning parametrisations $\phi_p$ point-by-point on a sphere by
starting at a particular point and spreading out in all directions.
The non-flatness of the sphere causes these initially divergent
directions to begin to converge, with some points ultimately reached
from multiple directions.  Unless special care is taken, the
parametrisations will not match up at such points.

This has not been seen before as a problem because implicit or
explicit use typically has been made of a Riemannian metric to guide
the construction of parametrisions.  Furthermore, the affine
invariance of the Newton method plays a critical role as it means
parametrisations constructed locally only have to agree globally
with each other up to affine transformations, which is easier
(indeed, possible) to achieve.  When thinking of a manifold as an
object embedded in Euclidean space, it is visually clear how the
affine tangent plane can be moved around the manifold, and although
two different paths from $p$ to $q$ may move the affine tangent
plane differently, the only difference will be a rotation.  Strategies
such as projection from the affine tangent plane onto the
manifold~\cite{Manton:opt_mfold}, or the use of the Riemannian
exponential map, can then be used to generate a smooth parametrisation
$\phi\colon TM \rightarrow M$ with which the Newton method can be
lifted from Euclidean space to the manifold $M$.

The sphere $S^2$ highlights the role of affine invariance.  Although
smooth parametrisations $\phi\colon TS^2 \rightarrow S^2$ are readily
constructed, the hairy-ball theorem states there is no global section
of the frame bundle on $S^2$, meaning there is no way to identify
each $T_pS^2$ with $\reals^2$ in a smooth way.  If not for the
affine invariance of the Newton method, it would be impossible to
construct smoothly varying parametrisations with which to lift the
Newton method to $S^2$.  For more general iteration functions,
this highlights the importance of conditions C1--C2 not requiring
continuity; see \refapp{cond}.

\section{Whimsical and Path-Dependent Newton Methods}
\label{sec:whim}

The global topological constraint in \refsec{gtc} makes it desirable
to allow Newton methods on manifolds to construct their parametrisations
dynamically as the sequence of iterates unfolds.  Such methods are
called path-dependent Newton methods because the parametrisation
$\phi_{p_k}$ used at the $k$th step may depend on the path
$p_0,p_1,\cdots,p_{k-1}$ leading up to $p_k$.  One of the many
possibilities this opens up is using (non-metric) affine connections
and parallel transport to construct parametrisations.  See
\refsec{recent} for other possibilities.

The theory in \refsec{nmm} extends to encompass path-dependent
Newton methods because there is no inherent requirement for the
$\phi_p$ to vary smoothly, or even continuously, in $p$.  The essence
of H1--H3 in \refsec{nmm} is that the bounds $\alpha$ and $\beta$
hold \textit{uniformly} on sufficiently small neighbourhoods.
Therefore, a generalised Newton method can be constructed by
specifying a $\phi_p$ for each $p$ with little regard for how the
$\phi_p$ fit together to form $\phi$.  In other words, a generalised
Newton method at each step is free to choose from many different
parametrisations $\phi_p$ without affecting its performance.

It is expedient to study path-dependent Newton methods in terms
of whimsical Newton methods.
Let $\Sigma = \{ \Sigma_p \mid p \in M\}$ be a collection of sets
$\Sigma_p$ of pairs $(\phi_p,\psi_p)$ of parametrisations $\phi_p,\psi_p\colon
T_pM \rightarrow M$.  A whimsical Newton method with respect to
$\Sigma$ is the general term given to any iterative scheme $p_{k+1}
= E^k_f(p_k)$ where $E^k_f(p)$ is a generalised Newton iteration
function \refeq{gni} using a pair of parametrisations $(\phi_p,\psi_p)$
belonging to $\Sigma_p$.  Indexing $E_f$ by $k$ means that even if
$p_{k+j} = p_k$ for some $j>0$, a different parametrisation pair
can be chosen from $\Sigma_{p_k}$ at the $k$th and $(k+j)$th steps.

Imposing the following uniformity constraints on the elements of
$\Sigma$ ensures any whimsical Newton method with respect to
$\Sigma$ converges locally quadratically to non-degenerate
critical points. Recall the definition of $A_p$ in \refsec{nmm}.

Conditions E0--E2 are satisfied if $\Sigma_p$ is
non-empty for all $p \in M$, and furthermore,
$\forall \pstar \in M$,
$\exists$ $\Ct^2$-chart $(U,\varphi)$ with $\varphi(\pstar)=0$,
$\exists \rho > 0$,
$\exists \alpha, \beta \in \reals$,
$\forall p \in U$,
$\forall (\phi_p,\psi_p) \in \Sigma_p$:
\begin{description}
\item[E0] $\phi_p(0_p)=\psi_p(0_p)=p$ and 
    $D\phi_p(0_p) = D\psi_p(0_p) = I$;
\item[E1] $\widehat\phi_p=\varphi \circ \phi_p \circ A_p^{-1}$
satisfies $\|D^2\widehat\phi_p(0)\| \leq \alpha$;
\item[E2] $\widehat\psi_p=\varphi \circ \psi_p \circ A_p^{-1}$
satisfies $\|\widehat\psi_p(y) - \varphi(p) - y\| \leq \beta \|y\|^2$
whenever $\|y\| < \rho$.
\end{description}

\begin{thm}
\label{th:whim}
Let $f\colon M \rightarrow \reals$ be a $\Ct^2$-smooth cost function on
a $\Ct^2$-smooth manifold $M$ satisfying \refeq{nmmf} at a
non-degenerate critical point $\pstar$.  If $\Sigma$ satisfies
E0--E2 then any whimsical Newton method $p_{k+1} = E^k_f(p_k)$ with
respect to $\Sigma$ converges locally quadratically to $\pstar$.
Uniform bounds exist for the rate of convergence (that is, the
constants $\kappa$ and $K$ in \refeq{cnv} and \refdef{conv} with
respect to a given local coordinate chart) that are independent of
how $E^k_f(p)$ selects which pair $(\phi_p,\psi_p) \in \Sigma_p$ to use
at the $k$th step.
\end{thm}

\refth{whim} can be proved by observing in the proof of \refth{cad}
that the rate of convergence is determined purely in terms of bounds
on the second-order behaviour of the parametrisations.  Provided
the bounds $\alpha$ and $\beta$ remain valid, the pairs of
parametrisations used become irrelevant.  Similarly,
it follows from the proofs of \refprop{Cinv} and \reflem{coch} that
if conditions E0--E2 are satisfied with respect to one chart, they
are satisfied with respect to any other.

A path-dependent Newton method differs from a whimsical Newton
method in that the rule for choosing the parametrisation pair to
use at each step may depend on previous iterates.

\begin{cor}
Let $f\colon M \rightarrow \reals$ be a $\Ct^2$-smooth cost function on
a $\Ct^2$-smooth manifold $M$ satisfying \refeq{nmmf} at a
non-degenerate critical point $\pstar$.  Let $p_{k+1} =
E^{(k,p_0,\cdots,p_{k-1})}_f(p_k)$ be a path-dependent Newton iterate
with respect to an indexed family $\Sigma = \{\Sigma_p \mid p \in
M\}$ of sets $\Sigma_p$ of parametrisation pairs: $p_{k+1} =
\psi_{(k)} \circ N_{f \circ \phi_{(k)}}(0_{p_k})$ where the rule
for choosing $(\psi_{(k)},\phi_{(k)}) \in \Sigma_{p_k}$ may depend
on past iterates $p_0,\cdots,p_{k-1}$ as well as on $k$ and $p_k$.
If $\Sigma$ satisfies E0--E2 then this path-dependent Newton iterate
converges locally quadratically to $\pstar$.
\end{cor}
\begin{proof}
Assume to the contrary the existence of $f$, $\pstar$ and a sequence
of initial points $\{p^{(i)}_0\}$ converging to $\pstar$
such that the path-dependent Newton iterate started at any $p^{(i)}_0$
does not converge quadratically to $\pstar$.  For each $p^{(i)}_0$,
the resulting path-dependent Newton iterate is a whimsical Newton
method with respect to $\Sigma$.  From \refth{whim}, there exists
a neighbourhood of $\pstar$ such that \textit{any} whimsical Newton
iterate with respect to $\Sigma$ that starts within this neighbourhood
will converge quadratically to $\pstar$, a contradiction.
\end{proof}

The motivation given earlier for introducing path-dependent methods
was that it is easier and more natural to construct parametrisations
locally then extend path-wise than to construct parametrisations
globally because the latter requires the local parametrisations to
fit together globally.  Another use for path-dependent methods is
to give the algorithm memory.  This leads into the study of general
techniques for extending conjugate gradient and other such methods to
manifolds, a topic outside the scope of the present paper.

\section{Re-Centring and Other Parametrisation Construction Techniques}
\label{sec:recent}

Various strategies exist for choosing parametrisations that satisfy
C1--C2 in \refsec{nmm}.  If the class of cost functions of interest
is known beforehand then this knowledge should inform the choice
of parametrisation; see \refsec{cad_disc}.  There is also interest
in choosing relatively simple parametrisations leading to generic
algorithms designed without regard to any particular class of cost
functions.  A basic idea for how to do this is introduced in
\refsec{sym} and generalised in subsequent sections.  It is called
re-centring and exploits the existence of a local diffeomorphism
between any two parts of a manifold.  \textit{This changes the focus
from devising parametrisations to devising transformations.} If $M$
were a sphere, for example, then instead of producing a sequence
of points $p_1,p_2,\cdots$ on $M$ converging to a critical point
$\pstar$, the manifold $M$, along with the cost function $f$, can be
rotated at each step to bring $p_k$ to the North pole, until
eventually the critical point $\pstar$ is brought to the North pole.
Since each Newton step is always taken from the North pole, its
design is simplified.

\subsection{Submersions and Fibre Bundles}
\label{sec:fb}

It may happen that the manifold $N$ for which a parametrisation is
sought is the image of a smooth function $g\colon M \rightarrow N$ where
$M$ is simpler to parametrise. For example, $M$ might be a matrix
Lie group and $N$ a homogeneous space.  Since a cost function $f\colon
N \rightarrow \reals$ pulls back to a cost function $f \circ g$ on
$M$, an iterative scheme on $M$ should induce an iterative scheme
on $N$. Simply pulling $f$ back is not recommended if $\dim M >
\dim N$ because the final algorithmic complexity might increase and
non-degenerate critical points of $f$ can become degenerate critical
points of $f \circ g$.

An alternative is to endeavour to ``push forwards'' the parametrisation
on $M$.  Let $g_\star\colon TM \rightarrow TN$ be the induced push-forward
of $g\colon M \rightarrow N$ and let $p \in M$ be such that $g$ is a
submersion at $p$.  Then $T_pM$ splits into a vertical component
$V_p$ and a non-unique horizontal component $H_p$, that is, there
exists a subspace $H_p \subset T_pM$ such that $T_pM = V_p \oplus
H_p$ where $V_p = \{v \in T_pM \mid g_\star(v)=0\}$.  Since $g$ is
a submersion at $p$, $g_\star$ induces a linear isomorphism from
$H_p$ to $T_{g(p)}N$, denoted $g_\star|_{H_p}$.  In particular, a
parametrisation $\phi_p\colon T_pM \rightarrow M$ can be used to form
the parametrisation $\tilde\phi_{g(p)} = g \circ \phi_p \circ
(g_\star|_{H_p})^{-1}$ from $T_{g(p)}N$ into $N$.

Take $\Sigma_q$ to be the set of all parametrisation pairs
$(\tilde\phi_q,\tilde\psi_q)$ coming from parametrisations
$\phi_p,\psi_p\colon T_pM \rightarrow M$ where $q = g(p)$ and $g$
is a submersion at $p$.  Provided every point $q \in N$ has at least
one preimage $p$ such that $g$ is a submersion at $p$, a whimsical
or path-dependent Newton method is well-defined; see \refsec{whim}.
Assuming $\phi$ and $\psi$ satisfy C1--C2 in \refsec{nmm}, it is
not necessarily the case though that the resulting $\Sigma$ will satisfy
E1--E2 in \refsec{whim}.  There are essentially two ways for E1--E2
to fail to hold.  Visually, the first is if the angle (with respect
to a Riemannian metric placed on $M$) between $H_p$ and $V_p$ can
approach zero as $p$ varies.  The second is if $g^{-1}(q)$ is
unbounded and $\phi_p$ for $p \in g^{-1}(q)$ gets arbitrarily
ill-behaved as $p$ goes to infinity.

If $g\colon M \rightarrow N$ happens to be a compact fibre bundle
and an Ehresmann connection is chosen, thereby determining the
horizontal bundle, then a compactness argument can be made
for E1--E2 to hold.

A more general approach is to limit the number of parametrisation
pairs in $\Sigma_q$ by judiciously choosing which preimages $p \in
g^{-1}(q)$ to use.  For example, if there exist a locally finite
open cover $\{U_\gamma\}$ of $N$ and smooth functions $h_\gamma\colon
U_\gamma \rightarrow M$ such that each $g \circ h_\gamma$ is the
identity map, then $\Sigma_q$ need only contain the finite number
of parametrisation pairs coming from the preimages $p = h_\gamma(q)$
for those $\gamma$ for which $q \in U_\gamma$, and C1--C2 will imply
E1--E2.

\subsection{Re-centring via a Group Action}
\label{sec:sym}

If a Lie group $G$ acts transitively on $M$ then a generalised
Newton method can be devised by continually re-centring the cost
function about a distinguished point.  Precisely, fix $\bar p \in
M$, choose parametrisations $\phi_{\bar p}, \psi_{\bar p}\colon T_{\bar
p}M \rightarrow M$ and define $E_f(\bar p)$ as in \refeq{gni}.  The
group action $p \mapsto g \cdot p$ allows $E_f(p)$ to be defined
for $p \neq \bar p$ by $E_f(p) = h(p) \cdot E_{\tilde
f}(\bar p)$ where $h\colon M \rightarrow G$ is an arbitrary
function satisfying $h(p) \cdot \bar p = p$, and
$\tilde f(q) = f(h(p) \cdot q)$ is the re-centred cost function.

If no obvious rule for choosing $h$ comes to mind, the framework
of \refsec{whim} can be used.  Precisely, define $\Sigma_p = \{
(\phi_g,\psi_g) \mid g \in G,\, g \cdot \bar p = p\}$ where $\phi_g
= \theta_g \circ \phi_{\bar p} \circ Q_g^{-1}$, $\psi_g = \theta_g
\circ \psi_{\bar p} \circ Q_g^{-1}$, $\theta_g(p) = g \cdot p$ and
$Q_g$ is the restriction of $(\theta_g)_\star$ to $T_{\bar p}M$.
Here, $(\theta_g)_\star\colon TM \rightarrow TM$ is the push-forward of
$\theta_g\colon M \rightarrow M$.  The parametrisation pair $(\phi_g,\psi_g)$
yield the same Newton step as before: $\psi_g \circ N_{f \circ
\phi_g}(0_{g \cdot \bar p}) = g \cdot E_{\tilde f}(\bar p)$ where
$\tilde f(q) = f(g \cdot q)$.  Assume $\phi_{\bar p}$ and $\psi_{\bar
p}$ are sensible (e.g., $\Ct^2$-smooth).  Provided the first and
second order derivatives of $\theta_g$ are bounded in an appropriate
sense then E0--E2 will hold; see \reflem{symh}.  This will be the
case if $G$ is compact, for example.  Otherwise, the number of
parametrisations in each $\Sigma_p$ can be limited with the aid of
a finite open cover of $M$ on which local sections are defined;
see the end of \refsec{fb}.

\subsection{Re-centring via Affine Transformations}
\label{sec:canon}

Re-centring can be applied to a manifold $M$ embedded in Euclidean
space by using affine transformations of Euclidean space to bring
any point $p \in M$ of interest to the origin in such a way that
the transformed version of $M$ is a graph of a function in a
neighbourhood of the origin.  Therefore, a rule for parametrising
a graph of a function induces a rule for parametrising $M$.  This
idea will be used in \refsec{proj}.  The present section focuses
on the reverse direction: determine if a parametrisation pair
$(\phi, \psi)$ satisfies C1--C2 in \refsec{nmm} by studying the
corresponding re-centred parametrisations.

\newcommand{\cH}{\mathfrak{H}}

Denote by $\cH_{\delta}$ the space of $\Ct^2$-smooth functions $h\colon
B_n(0;\delta) \rightarrow \reals^k$ satisfying $h(0)=0$, $Dh(0)=0$
and $\sup_{t \in B_n(0;\delta)} \| D^2h(t) \| < \infty$.  (Recall
$B_n$ is an open ball in $\reals^n$.) Associate with any $h \in
\cH_\delta$ the manifold
\begin{equation}
\label{eq:Mh}
\tilde M_{(h,\delta)} = \{(t,h(t)) \mid t \in B_n(0;\delta)\} \subset
\reals^{n+k}.
\end{equation}

A parametrisation $\tilde\phi_{(0,0)}\colon T_{(0,0)}\tilde M_{(h,\delta)}
\rightarrow \tilde M_{(h,\delta)}$
can be represented by a function $\pi\colon \reals^n \rightarrow
\reals^n$ taking $x$ to the point $(\pi(x),h \circ \pi(x))$ on
$\tilde M_{(h,\delta)}$.  Here, $x$ represents the point $(x,0)$
on the affine tangent space at the origin of $\tilde M_{(h,\delta)}$.
An exemplar is using Euclidean projection from $(x,0)$ to $\tilde
M_{(h,\delta)}$ to define $\pi$, so that $\pi$ satisfies
\begin{equation}
\label{eq:expirule}
\| (\pi(x), h \circ \pi(x)) - (x,0) \| = \min_{t \in B(0;\delta)}
\| (t, h(t)) - (x,0) \|
\end{equation}
whenever the minimum exists.  This will be studied in \refsec{proj}.

It is not important for $\pi$ to be defined uniquely by a rule, nor
for $\pi$ to be defined on the whole of $\reals^n$.  Essentially,
it is merely required that any choice of $\pi$ is defined on a
sufficiently small domain $B_n(0;\rho)$ and satisfies $\| D^2\pi(0) \| \leq
\alpha$ and $\| \pi(x) - x \| \leq \beta \|x\|^2$ for $x \in
B_n(0;\rho)$, where the constants $\rho$, $\alpha$ and $\beta$
depend on $h$ and $\delta$ in a way that ensures they remain uniformly
bounded if $h$ and $\delta$ are perturbed; see \refsec{whim}.
This is now made precise.

A parametrisation $\phi\colon TM \rightarrow M$ of an $n$-dimensional
embedded submanifold $M \subset \reals^{n+k}$ is said to satisfy
condition \textbf{D} if D1--D6 below are satisfied.  Central to
this condition is a class of functions $\pi$ obtained from $\phi$
as follows.  Associate to each $x \in M \subset \reals^{n+k}$ a
rotation $R_x$ of $\reals^{n+k}$ sending $V_xM$ to $\reals^n \times
\{0\}$.  (Here, the notation of \refsec{embed} is being used, but
with $i$ and $i_\star$ omitted because $M$ is being treated as an
actual subset of $\reals^{n+k}$.) Define the translation $Q_x(q)=q-x$.
For any pair $(h,\delta)$ for which $R_xQ_x(M)$ locally looks like
$\tilde M_{(h,\delta)}$, meaning there exists a set $U$ open in
$\reals^{n+k}$ such that $R_xQ_x(M) \cap U = \tilde M_{(h,\delta)}$,
a function $\pi\colon \reals^n \rightarrow \reals^n$ can be defined
implicitly by
\begin{equation}
\label{eq:pi}
(\pi(y), h \circ \pi(y)) = R_x\, Q_x\, \widehat\phi(x,R_x^{-1}Jy)
\end{equation}
where $\widehat\phi = \phi \circ \tau_I^{-1}$ is the representation
of $\phi$ in local coordinates with respect to the identity chart
(on $\reals^{n+k}$; refer to \refsec{embed}) and $J\colon \reals^n
\rightarrow \reals^{n+k}$ sends $y$ to $J(y)=(y,0)$.  To emphasise,
$\pi$ is only defined at points $y \in \reals^n$ for which \refeq{pi}
holds.  Note too that $\pi$, which depends on the triple $x$, $h$
and $\delta$, is merely $\phi_x$, the restriction of $\phi$ to
$T_xM$, written in a canonical form (albeit depending on the choice
of rotation $R_x$).

Condition D requires there to exist functions $\alpha, \beta,
\rho\colon (0,\infty) \times [0,\infty) \rightarrow \reals$ such that
for all $x \in M$, $\delta > 0$ and $h \in \cH_\delta$ for which
$R_xQ_x(M)$ locally looks like $\tilde M_{(h,\delta)}$ and
for $K = \sup_{t \in B_n(0;\delta)} \|D^2h(t)\|$:
\begin{description}
\item[D1] The domain of definition of $\pi$ in \refeq{pi} includes 
$B_n(0;\rho(\delta,K))$;
\item[D2] $\| D^2\pi(0) \| \leq \alpha(\delta,K)$;
\item[D3] $\| \pi(y) - y \| \leq \beta(\delta,K)\, \|y\|^2$ for
$y \in B_n(0;\rho(\delta,K))$.
\end{description}

It is also required that for all $\delta > 0$ and all
$\bar K \in [0,\infty)$:
\begin{description}
\item[D4] $\sup_{K \in [0,\bar K]} \alpha(\delta,K) < \infty$;
\item[D5] $\sup_{K \in [0,\bar K]} \beta(\delta,K) < \infty$;
\item[D6] $\inf_{K \in [0,\bar K]} \rho(\delta,K) > 0$.
\end{description}

A sufficient condition for D4--6 to hold is for $\alpha$ and $\beta$
to be upper semi-continuous in $K$, and $\rho$ lower semi-continuous
in $K$.

\begin{prop}
\label{prop:canon}
Let $\phi\colon TM \rightarrow M$ be a parametrisation of an $n$-dimensional
$\Ct^2$-smooth embedded submanifold $M \subset \reals^{n+k}$.  If $\phi$
satisfies condition D described above then $\phi$ and $\psi=\phi$
satisfy C1--C2 in \refsec{nmm}.  (If $\phi$ satisfies D1, D2, D4
and, instead of D3 and D6, the weaker conditions that $\pi(0)=0$,
$D\pi(0)=I$ and $\rho(\delta,K) > 0$, then $\phi$ satisfies C1.  If
$\phi$ satisfies D1, D3, D5 and D6 then $\psi=\phi$ satisfies C2.)
\end{prop}
\begin{proof}
When convenient, elements of $\reals^{n+k}$ are written as $(u,v)
\in \reals^n \times \reals^k$, with projections $P_1$ and $P_2$
sending $(u,v)$ to $u$ and $v$ respectively.  Fix a point $z \in M
\subset \reals^{n+k}$.  Let $\delta > 0$ and $h \in \cH_{6\delta}$
be such that $\tilde M_{(h,6\delta)} = R_zQ_z(M) \cap U$ for some
open set $U \subset \reals^{n+k}$.  Let $K = \sup_{t \in B(0;6\delta)}
\|D^2h(t)\|$.  By shrinking $\delta$ if necessary, it is assumed
without loss of generality that $\delta < 1/(4K)$ and $\| Dh(t) \|
< 1/4$ for $t \in B(0;6\delta)$.

Choose an arbitrary $t \in B(0;\delta)$.  Define $x = Q_z^{-1}
R_z^{-1}(t,h(t))$ and $f_t(\tau) = P_1 R_x R_z^{-1}(\tau,
h(t+\tau)-h(t))$ for $\tau \in B(-t;6\delta) \supset B(0;5\delta)$.
It will be shown that $R_xQ_x(M)$ locally looks like $\tilde
M_{(h_t,\delta)}$ where $h_t(u) = P_2 R_x R_z^{-1}(f_t^{-1}(u),
h(t+f_t^{-1}(u)) - h(t))$.

First, bounds on $f_t$ and its derivatives $Df_t(\tau) \cdot \xi =
P_1 R_x R_z^{-1} (\xi, Dh(t+\tau) \cdot \xi)$ and $D^2f_t(\tau)
\cdot \xi^2 = P_1 R_x R_z^{-1} (0, D^2h(t+\tau) \cdot \xi^2)$ are
obtained.  Importantly, $P_2 R_x R_z^{-1} (\xi, Dh(t) \cdot \xi) =
0$ because $R_z^{-1} (\xi, Dh(t) \cdot \xi)$ lies in $V_x M$, the
latter a consequence of the curve $\gamma(\theta) = Q_z^{-1} R_z^{-1}
(t+\theta\xi, h(t+\theta\xi))$ lying on $M$ and having $\gamma'(0)
= R_z^{-1}(\xi, Dh(t) \cdot \xi)$.  Therefore, $\| Df_t(0) \cdot
\xi \| = \| (\xi, Dh(t) \cdot \xi) \| \geq \| \xi \|$.  Similarly,
$\| f_t(\tau) \| \geq \| \tau \| - \| h(t+\tau) - h(t) - Dh(t)\cdot\tau
\| \geq \| \tau \| - \frac12 K \| \tau \|^2$ (see \reflem{lip}).  For
$\tau_1,\tau_2 \in B(-t;6\delta)$, $\| Df_t(\tau_1) - Df_t(\tau_2)
\| \leq \| Dh(t+\tau_1) - Dh(t+\tau_2) \| < \frac12$.  Thus $\|
(Df_t(\tau))^{-1} \| < 2$ because, for $\| \xi \| > 0$, $\|Df_t(\tau)
\cdot \xi\| \geq \|Df_t(0) \cdot \xi\| - \|(Df_t(\tau) - Df_t(0))
\cdot \xi\| > \frac12 \|\xi\|$.  Lemma 1 of~\cite[Chapter 16]{bk:Hirsch:diff}
implies $f_t$ is injective on $B(-t;6\delta)$, and, since $Df_t$
is invertible, the inverse function theorem implies $f_t$ is a
$\Ct^2$-diffeomorphism from $B(-t;6\delta)$ onto its image.
Furthermore, since $\| \tau \| = 4 \delta$ implies $\| f_t(\tau) -
f_t(0) \| \geq 4\delta - \frac12 K(4\delta)^2 > 2 \delta$, the proof
of Lemma 2 of~\cite[Chapter 16]{bk:Hirsch:diff} implies $B(0;\delta)
\subset f_t(B(0;4\delta))$.  Thus, $h_t(u)$ is a well-defined
$\Ct^2$-smooth function on $B(0;\delta)$.  Finally, note $\|
D^2f_t(\tau) \| \leq \| D^2h(t+\tau) \| \leq K$.

Next it is shown that $\tilde M_{(h_t,\delta)} = R_x R_z^{-1}
Q_{(t,h(t))} (\tilde M_{(h,6\delta)}) \cap P_1^{-1}( B(0;\delta)
)$.  Indeed, if $(u,h_t(u)) \in \tilde M_{(h_t,\delta)}$ then $\tau
= f_t^{-1}(u)$ is well-defined and $(t+\tau,h(t+\tau)) \in \tilde
M_{(h,6\delta)}$ is such that $(u,h_t(u)) = R_x R_z^{-1} Q_{(t,h(t))}
(t+\tau,h(t+\tau))$.  Conversely, an arbitrary element $p$ of $R_x
R_z^{-1} Q_{(t,h(t))} (\tilde M_{(h,6\delta)})$ is of the form $p
= R_x R_z^{-1}(\tau, h(t+\tau)-h(t))$ for some $\tau \in B(-t;6\delta)$.
Since $P_1(p) = f_t(\tau)$, if $p \in P_1^{-1}(B(0;\delta))$ then
$u = f_t(\tau) \in B(0;\delta)$ and $p = (u,h_t(u)) \in \tilde
M_{(h_t,\delta)}$.  Because $(t,h(t)) = R_z Q_z (x)$, it follows
that $R_x R_z^{-1} Q_{(t,h(t))} R_z Q_z = R_x Q_x$.  Therefore,
$\tilde M_{(h_t,\delta)} = R_x Q_x(M) \cap R_x Q_x Q_z^{-1} R_z^{-1}(U)
\cap P_1^{-1}(B(0;\delta))$, proving $R_x Q_x(M)$ locally looks
like $\tilde M_{(h_t,\delta)}$.

To show $h_t \in \cH_{\delta}$, first note $h_t(0) = 0$ and $Dh_t(0)
= 0$, the latter a consequence of $P_2 R_x R_z^{-1} (\xi, Dh(t)
\cdot \xi) = 0$ for all $\xi$. It remains to bound $\| D^2h_t(u)
\|$.  For an arbitrary $u \in B(0;\delta)$, let $\tau = f_t^{-1}(u)
\in B(0;4\delta)$.  Then $D^2h_t(u) \cdot \xi^2 = P_2 R_x R_z^{-1}
(D^2f_t^{-1}(u) \cdot \xi^2, D^2h(t+\tau) \cdot (Df_t^{-1}(u) \cdot
\xi)^2 + Dh(t+\tau) \, D^2f_t^{-1}(u) \cdot \xi^2)$.  Now, $Df_t^{-1}(u)
= (Df_t(\tau))^{-1}$ and $D^2 f_t^{-1}(u) \cdot \xi^2 = -
(Df_t(\tau))^{-1} \, D^2f_t(\tau) \cdot ( (Df_t(\tau))^{-1} \cdot
\xi)^2$.  Applying the earlier bounds shows $\| D^2h_t(u) \| \leq
\bar K$ where $\bar K = \sqrt{ (8K)^2 + (4K + 2K)^2 } = 10K$.

With the values of $\delta$ and $\bar K$ as above, and referring
to D4--D6, define $\bar\alpha = \sup_{K \in [0,\bar K]} \alpha(\delta,K)$,
$\bar\beta = \sup_{K \in [0,\bar K]} \beta(\delta,K)$, and $\bar\rho
= \inf_{K \in [0,\bar K]} \rho(\delta,K)$.  Shrink $\bar\rho$ if
necessary to ensure $\bar\rho < \delta$ and $B_{n+k}(0;2\bar\rho)
\subset U$.  Thus, for any $x \in M \cap B_{n+k}(z;\bar\rho)$ there
exists a $t \in B(0;6\delta)$ such that $x = Q_z^{-1} R_z^{-1} (t,
h(t))$, in which case $\|t\| \leq \|R_z^{-1} (t,h(t))\| = \| x - z
\| < \bar\rho < \delta$.  In particular then, $R_xQ_x(M)$ locally
looks like $\tilde M_{(h_t,\delta)}$.  For the particular triple
$x$, $h_t$ and $\delta$, define $\pi$ as in \refeq{pi}.  From D1--D3
it follows that $\pi$ is defined on $B(0;\bar\rho)$ and satisfies
$\|D^2\pi(0)\| \leq \bar\alpha$ and $\| \pi(\tilde y) - \tilde y
\| \leq \bar\beta\,\|\tilde y\|^2$ for $\tilde y \in B(0;\bar\rho)$.

The proof is completed by showing $\widehat\phi$ satisfies the
assumptions in \reflem{affine_embed}.  Let $y \in B(0;\bar\rho)
\cap V_x M$ be arbitrary.  Define $\tilde y = P_1 R_x y$ and note
$R_x y = (\tilde y, 0)$.  From D1 and \refeq{pi}, since $\| \tilde
y \| = \| y \| < \bar\rho$, $\widehat\phi_x(y) = Q_x^{-1} R_x^{-1}
(\pi(\tilde y), h_t \circ \pi(\tilde y))$.  Then $\widehat\phi_x(y)$
lies in $M$ because $R_x Q_x(M)$ locally looks like $\tilde
M_{(h_t,\delta)}$.  That $\widehat\phi_x(0)=0$ follows from $\pi(0)=0$
(D3) and $h_t(0)=0$.  The facts $\pi(0)=0$, $D\pi(0)=I$ (D3) and
$Dh_t(0)=0$ imply $D\widehat\phi_x(0) \cdot y = y$ and $D^2\widehat\phi_x(0)
\cdot y^2 = R_x^{-1} ( D^2\pi(0) \cdot \tilde y^2, D^2h_t(0) \cdot
\tilde y^2 )$.  Thus $\| D^2 \widehat\phi_x(0) \cdot y^2 \| \leq
\sqrt{\bar\alpha^2 + \bar K^2} \, \|y\|^2$.  \reflem{lip} implies
$\| h_t(u) \| \leq \frac12 \bar K \| u \|^2$.  Therefore, $\| h_t
\circ \pi(\tilde y) \| \leq \frac12 \bar K (\| \pi(\tilde y) - \tilde
y \|^2 + \| \tilde y\|^2 ) \leq \frac12 \bar K (1+\bar\beta^2\bar\rho^2)
\| \tilde y \|^2$.  Thus, $\| \widehat\phi_x(y) - x - y \| = \|
(\pi(\tilde y) - \tilde y, h_t \circ \pi(\tilde y)) \| \leq \sqrt{
\bar\beta^2 + (1/4) \bar K^2 (1+\bar\beta^2\bar\rho^2)^2 } \| y
\|^2$.
\end{proof}

\subsection{Local and Global Projections}
\label{sec:proj}

For an embedded manifold $M \subset \reals^{n+k}$, it was suggested
in~\cite{Manton:opt_mfold} that projection from the affine tangent
plane to the manifold could be used to define parametrisations.
The utility of \refprop{canon} is illustrated by proving such
parametrisations satisfy C1--C2.  Only ordinary calculus is required
as the differential geometric framework is hidden behind \refprop{canon}.
Note the proof works at the generality of $\Ct^2$-smooth manifolds
and is thus not based on a smoothness argument (\reflem{Mthree}).

\begin{lem}
\label{lem:hproj}
Let $h\colon B_n(0;\delta) \rightarrow \reals^k$ be a $\Ct^2$-smooth map
with $h(0)=0$ and $Dh(0)=0$.  Let $K = \sup_{x \in B(0;\delta)}
\|D^2h(x)\|$ and $\rho = \frac12 \min\{\delta, \sqrt{2}/(3K)\}$.
There exists a $\Ct^1$-smooth map $\pi\colon B_n(0;\rho) \rightarrow
B_n(0;2\rho)$ such that, for any $x \in B(0;\rho)$, $(\pi(x),h\circ\pi(x))$
is the unique point on the manifold $\{(t,h(t)) \mid t \in
B(0;\delta)\}$ closest to $(x,0)$. Moreover, no point is closest
to $(x,0)$ on the smaller manifold $\{(t,h(t)) \mid \|t\| <
\|\pi(x)\|\}$.  The map $\pi$ satisfies $\pi(0)=0$, $D\pi(0)=I$,
$D^2\pi(0)=0$ and $\|\pi(x)-x\| \leq \frac12 K \|x\|^2$.
\end{lem}
\begin{proof}
Define $f(x) = x + g(x)$ and $g(x) = (Dh(x))^T h(x)$ where superscript
$T$ denotes adjoint.  Then $(f(x),0)$ is the unique point of
intersection of the affine plane normal to the manifold at $(x,h(x))$
with the plane $\reals^n \times 0$.  Thus if $\pi$ exists it
must satisfy $f(\pi(x))=x$.

From \reflem{lip}, $\|h(x)\| \leq (K/2) \|x\|^2$ and $\|Dh(x)\|
\leq K \|x\|$ for $x \in B(0;\delta)$.  Therefore, $\| g(x) \| \leq
(1 / 2) K^2 \|x\|^3$, $\| Dg(x) \| \leq (3/2) K^2 \|x\|^2$, $\|
Df(x) \cdot \xi \| \geq (1 - (3/2) K^2 \|x\|^2) \|\xi\|$ and $\|
Df(x) - Df(y) \| \leq (3/2) K^2(\|x\|^2 + \|y\|^2)$ for $x,y \in
B(0;\delta)$.  The latter implies $D^2f(0)=0$.  If $x,y \in B(0;2\rho)$
then $\|(Df(x))^{-1}\| < 3/2$ and $\| Df(x) - Df(y) \| < 2/3$ so
Lemma 1 of~\cite[Chapter 16]{bk:Hirsch:diff} and the inverse function
theorem imply $f$ restricted to $B(0;2\rho)$ is a $\Ct^1$-diffeomorphism.

Let $x \in B(0;\rho)$.  Since the manifold includes the origin, a
distance $\|x\|$ away from $(x,0)$, the closest point(s) to $(x,0)$
on the original manifold are the same as the closest point(s) on
the smaller manifold $\{(t,h(t)) \mid t \in \overline{B}(0;2\|x\|)\}$.
The latter manifold is compact and hence a closest point exists.
Uniqueness follows from $f$ being injective on $B(0;2\rho) \supset
\overline{B}(0;2\|x\|)$; any closest point $(t,h(t))$ must satisfy
$f(t)=x$.  Since this is a local condition, it also means no point
is closest to $(x,0)$ on $\{(t,h(t)) \mid \|t\| < \|\pi(x)\|\}$.
Note that $\|\pi(x)\| \leq 2 \|x\|$.

The geometric bound $\|\pi(x)-x\| \leq \|h(x)\|$ implies $\|
\pi(x) - x\| \leq (K/2) \|x\|^2$, so $\pi(0)=0$ and $D\pi(0)=I$.
As $f$ is a $\Ct^1$-diffeomorphism on $B(0;2\rho)$,
$\pi(x) = f^{-1}(x)$ for $x \in B(0;\rho)$ is $\Ct^1$-smooth and
$D\pi = (Df \circ \pi)^{-1}$. For $\epsilon > 0$, choose $\delta >
0$ such that $\|y\| < 2\delta$ implies $\| y \| \leq 2 \| f(y) \|$
and $\| I - Df(y) \| < \epsilon \|y\| < \frac12$.  (This is
possible because $D^2f(0)=0$.) Then for $\|x\| < \min\{\rho,\delta\}$
and $y = \pi(x)$, $\| D\pi(x) - I \| = \| (I-(I-Df(y)))^{-1} - I
\| \leq (1-\|I-Df(y)\|)^{-1} - 1 < 2 \epsilon \|y\| \leq 4 \epsilon
\| f(y) \| = 4 \epsilon \| x \|$, proving $D^2\pi(0)=0$.
\end{proof}

It is mentioned tangentially that the cubic bound $\| \pi(x) - x
\| = \| f(\pi(x)) - \pi(x) \| \leq 4 K^2 \|x\|^3$ is readily
obtainable from the above proof.

Given a rule such as \refeq{expirule}, a parametrisation of a
manifold can be obtained by using \refeq{pi} in reverse.  The only
technicality is the choice of neighbourhood size $\delta$ to use
for each point on the manifold.  The actual choice is generally not
important provided a positive uniform lower bound exists on any
compact neighbourhood.  In fact, as is the case in the following
lemma, the choice may depend on $y$ as well as $x$ in \refeq{pi}.

\begin{prop}
\label{prop:localproj}
Let $M \subset \reals^{n+k}$ be an embedded $\Ct^2$-manifold of
dimension $n$.  Adopting the notation of \refsec{canon}, associate
with each $x \in M$ and $y \in \reals^n$ a $\delta_{xy} > 0$, an
$h_{xy} \in \cH_{\delta_{xy}}$ and a rotation $R_x$ such that: 1)
$R_xQ_x(M)$ locally looks like $\tilde M_{(h_{xy},\delta_{xy})}$;
and, 2) if no point of $\tilde M_{(h_{xy},\delta_{xy})}$ is closest
to $(y,0)$ then the same is true for any admissible choice of
$\delta_{xy}$.  Referring to \refeq{pi} and \refeq{expirule}, if
$(t,h_{xy}(t)) \in \tilde M_{(h_{xy},\delta_{xy})}$ is the unique
closest point to $(y,0)$ then set $\widehat\phi(x,R_x^{-1}Jy) =
Q_x^{-1} R_x^{-1}\, (t,h_{xy}(t))$.  Otherwise, if the closest point
does not exist or is not unique, let $\widehat\phi(x,R_x^{-1}Jy)$
be an arbitrary element of $M$.  Then the parametrisation $\phi\colon
TM \rightarrow M$, $\phi = \widehat\phi \circ \tau_I$, satisfies
condition D.
\end{prop}
\begin{proof}
The $k=0$ case is straightforward so assume $k>0$.  Define the
functions $\rho(\delta,K) = \frac12 \min\{\delta, \sqrt{2}/(3K)\}$,
$\alpha(\delta,K) = 0$ and $\beta(\delta,K) = K/2$; they satisfy
D4--D6.  Next, choose $x \in M$, $\delta > 0$ and $h \in \cH_{\delta}$
such that $R_x Q_x(M)$ locally looks like $\tilde M_{(h,\delta)}$.
Let $K = \sup_{t \in B(0;\delta)} \|D^2h(t)\|$.  Then for $y \in
B(0; \rho(\delta,K))$ define $\pi(y)$ as in \reflem{hproj}; the
unique closest point to $(y,0)$ on $\tilde M_{(h,\delta)}$ is
$(\pi(y),h\circ\pi(y))$.  This must therefore correspond with the
$\pi$ in \refeq{pi}.  That D1--D3 hold follows immediately from
\reflem{hproj}.
\end{proof}

To assist in interpreting \refprop{localproj}, consider how the
projection from $(y,0)$ onto $\tilde M_{(h,\delta)}$ changes with
$\delta$.  If $\delta$ is too small then no point is necessarily
closest because $\tilde M_{(h,\delta)}$ is not compact.  As $\delta$
increases the closest point may change as more candidates become
available.  The advantage of \refprop{localproj} in practice is it
allows parametrisations to be defined using only local minima of
the Euclidean distance function rather than insisting on global
minima.  Note too that projecting onto $\tilde M_{(h,\delta)}$ is
different from projecting onto $M$ because $M$ may curve around and
come close to touching itself.

\begin{prop}
\label{prop:proj}
Let $M \subset \reals^{n+k}$ be an embedded $n$-dimensional
$\Ct^2$-manifold.  Let $a\colon TM \rightarrow \reals^{n+k}$ be the map
taking a tangent vector to its equivalent point on the affine tangent
plane.  Let $\phi\colon TM \rightarrow M$ be any map with the property
that $\| \phi(v_p) - a(v_p) \| = \min_{q \in M} \| q - a(v_p) \|$
whenever the minimum exists, where the norm is the Euclidean norm
on $\reals^{n+k}$.  Then $\phi$ and $\psi=\phi$ satisfy C1 and C2
of \refsec{nmm}.
\end{prop}
\begin{proof}
The $k=0$ case is straightforward so assume $k>0$.  For $z \in M$,
define $\rho(z) = \sup\{ \rho \mid \exists h \in \cH_\rho,  \exists
\ \textrm{an open}\ U \supset B_{n+k}(0;\rho), R_zQ_z(M) \cap U =
\tilde M_{(h,\rho)}\}$; see \refsec{canon} for notation.  Let $K
\subset M$ be a compact set and assume to the contrary there
exists a convergent sequence $z_i \rightarrow \bar z$ in $K$ with
$\rho(z_i) \rightarrow 0$.  It follows from the proof of \refprop{canon}
that at $\bar z$ there exist a $\delta > 0$ and an $h \in \cH_{6\delta}$
such that (by shrinking $\delta$ if necessary) $\tilde M_{(h,6\delta)}
= R_{\bar z} Q_{\bar z}(M) \cap U$, $U \supset B_{n+k}(0;6\delta)$
and for any $z \in B(\bar z; \delta) \cap M$, there exists an $h_z
\in \cH_\delta$ such that $\tilde M_{(h_z,\delta)} = R_z Q_z(M)
\cap R_z Q_z Q_{\bar z}^{-1} R_{\bar z}^{-1}(U) \cap
P_1^{-1}(B_n(0;\delta))$.  Since $B_{n+k}(0;\delta) \subset R_z Q_z
Q_{\bar z}^{-1} R_{\bar z}^{-1}( B_{n+k}(0;6\delta))$ it follows
that $\rho(z) > \rho$ for $\|z - \bar z\| < \rho$, a contradiction.
Thus $\inf_{z \in K} \rho(z) > 0$.

For $x \in M$ and $y \in B_n(0;\rho(x)/4)$ define $\delta_{xy} =
\rho(x) / 2$ and $\widehat\phi(x,R_x^{-1}Jy) = \phi \circ
\tau_I^{-1}(x,R_x^{-1}Jy)$.  The closest point to $(y,0)$ on
$R_xQ_x(M)$ must be contained in $B_{n+k}(0;\rho(x)/2)$ and hence
is in the local representation $\tilde M_{(h_{xy},\delta_{xy})}$.
It is therefore possible to define $\delta_{xy}$ and
$\widehat\phi(x,R_x^{-1}Jy)$ for $y \not\in B_n(0;\rho(x)/4)$ so
that $\widehat\phi$ satisfies the conditions of \refprop{localproj}
(and hence by \refprop{canon} the corresponding parametrisation
satisfies C1 and C2).  By \reflem{rp}, because $\phi(v_x) =
\widehat\phi \circ \tau_I(v_x)$ whenever $\|v_x\| < \rho(x)/4$, 
$\phi$ and $\psi$ satisfy C1 and C2.
\end{proof}

Since any manifold can be embedded in $\reals^{n+k}$ for sufficiently
large $k$, \refprop{proj} guarantees the existence of parametrisations
satisfying C1--C2.

\subsection{Discussion}
\label{sec:gnm_disc}

To the best of our knowledge, all Newton methods on
finite-dimensional manifolds in the literature can be rewritten as
\refeq{gni} where the parametrisations $\phi$ and $\psi$ are smooth.
\refth{nmm} and \reflem{Mthree} together imply that such Newton
methods have local quadratic convergence.  As a specific example,
the original Riemannian Newton method in~\cite{Gabay:optm} uses the
Riemannian exponential map for the parametrisations $\phi$ and
$\psi$; see \refeq{riem_n}.  It is a standard result that if $M$
is $\Ct^4$-smooth then $\Exp$ is $\Ct^2$-smooth on a neighbourhood
of the zero section, and moreover, $\Exp_p(0)=p$ and $D\Exp_p(0)=I$.
Therefore, the Riemannian Newton method \refeq{riem_n} has local
quadratic convergence by \refrem{Mthree}.

The article~\cite{Manton:opt_mfold} introduced Newton methods on
(real and complex) Grassmann and Stiefel manifolds, with parametrisations
chosen to be global projections from Euclidean space onto the affine
tangent planes of the Stiefel manifold, and an analogous choice
made for the Grassmann manifold by treating it as a quotient space
of the Stiefel manifold.  Local quadratic convergence follows from
\refprop{proj} and \refsec{fb}.

When sufficient smoothness is not present for \reflem{Mthree} to
be applicable, the proofs in \refsec{recent} demonstrate that
essentially all the effort goes into obtaining \textit{uniform}
bounds.  Condition D in \refsec{canon} is one illustration of this.

The conjecture made in \refsec{canewt} applies to Newton methods
on manifolds too.  It is difficult to see how any iterative scheme can
fail to be of the form \refeq{gni} if it uses only the information
in the 2-jet of $f$ about the current point to converge locally
quadratically to a non-degenerate critical point for a sufficiently
rich class of functions $f$.

An advantage of expressing an algorithm in the form \refeq{gni} is
that it gives the algorithm the following heuristic interpretation:
at each step, the parametrisation $\phi$ endeavours to make $f \circ
\phi$ look as quadratic as possible, while $\psi$ endeavours to map
the result back to the manifold as cheaply as possible; see 
\refsec{cad_disc}.  Additionally, the fundamental idea of re-centring
can further simplify matters.

Finally, it is remarked that affine connections and parallel transport
can be used to construct parametrisations.  This can be understood
in terms of the classical notion of development in differential
geometry; a manifold $M$ can be rolled along an affine space without
slipping. This is a representative example of the re-centring
technique in \refsec{canon}.


\section{Iterates Computing Coordinate Independent Properties}
\label{sec:it}

This section studies how iteration functions besides the Newton
iteration function can be lifted from Euclidean space to manifolds.
This necessitates introducing a rudimentary theory of iterative
methods computing coordinate independent properties.
It also studies further the generalised Newton method at
a grass-roots level.

First, the concept of converging to an identifiable point of $f$
needs defining.  \refex{cip} may prove illuminative.

\newcommand{\Int}{\operatorname{Int}}
\newcommand{\Dom}{\operatorname{Dom}}

\begin{defn} \label{def:cind}
Assign to each $f\colon \operatorname{Dom} f \subset \rn
\rightarrow \reals$ a subset $P_f \subset \Int\Dom f$ of the
interior of the
domain of $f$.  The \emph{property} $P =
\{P_f\}$ is \emph{$\Ct^k$-coordinate independent} if $x \in P_f$
implies $\phi^{-1}(x) \in P_{f \circ \phi}$ for every $\Ct^k$-diffeomorphism
$\phi$ of open sets in $\reals^n$ with $x$ in the image of $\phi$.
\end{defn}

Henceforth, $f\colon \rn \rightarrow \reals$ will mean $f\colon
\operatorname{Dom} f \subset \rn \rightarrow \reals$ with
an implicit requirement that a particular point be in the domain
of $f$ whenever necessary.  For example, $f$ being $\Ct^k$-smooth
at $x$ implicitly requires $x \in \operatorname{Dom} f$.

Two functions $f,g\colon \rn \rightarrow \reals$ are $k$-jet equivalent
at $p \in \rn$ if $f$ and $g$ are $\Ct^k$-smooth in a neighbourhood
of $p$ and $f(p)=g(p)$, $Df(p)=Dg(p)$, $\cdots$, $D^kf(p) = D^kg(p)$.

\begin{defn} \label{def:fit}
A \emph{$k$th-order iterative method} is the assignment of an
iteration function $N_f\colon \rn \rightarrow \rn$ to each $f\colon \rn
\rightarrow \reals$ where $N_f(p) = N_g(p)$ whenever $f$ and $g$
are $k$-jet equivalent at $p$.
An iterative method \emph{computes the
property $P = \{P_f\}$ with
rate $K$} if, for any given $f$ and $\xstar \in P_f$, the iterate
$N_f$ converges locally to $\xstar$ with rate $K$.
\end{defn}

\begin{exmp} \label{ex:cip}
Let $P_f$ be the set of points
$x$ such that $f$ is $\Ct^3$-smooth in a neighbourhood of $x$, and
$x$ is a non-degenerate critical point of $f$.  Then $P$ is
$\Ct^3$-coordinate independent and the Newton iterate $N_f$ in \refeq{newt}
is a 2nd-order iterative method that computes $P$ with rate $2$.
\end{exmp}

\begin{exmp} \label{ex:cip2}
In \refex{cip}, $P_f$ can be instead the set of
points $\xstar$ satisfying the conditions in \refth{Neuclid},
including \refeq{Neuclid}.  
Then $P$ is a $\Ct^2$-coordinate independent property
(see \refrem{Neuclidind})
computed by the Newton iterate \refeq{newt}.
\end{exmp}

\begin{rem}
It follows from \refdef{cind} by using the identity map $\phi\colon U
\rightarrow U$ that if $P$ is a $\Ct^k$-coordinate independent
property, $U$ is an open subset of $\rn$ and $g=f|_U$ is the
restriction to $U$ of a function $f\colon \rn \rightarrow \reals$ then
$P_f \cap U \subset P_g$.  The converse is not implied; properties
can be ``forgotten'' as the domain increases.  Stricter definitions
precluding this are not necessary for what follows.
\end{rem}

Properties of functions in $\rn$ lift to properties of
functions on manifolds.

\begin{defn}
\label{def:Pbar}
Let $M$ be a manifold with maximal atlas $\mathcal{A}$ of
$\Ct^k$-smooth charts $(U,\varphi)$. Let $P = \{P_f\}$ be a
$\Ct^{k'}$-coordinate independent property with $k' \leq k$.  For
any $f\colon M \rightarrow \reals$, define $\bar P_f = \{p \in M \mid
\exists (U,\varphi) \in \mathcal{A},\,p\in U,\,\varphi(p) \in P_{f \circ
\varphi^{-1}}\}$.  The elements of $\bar P_f$ are said to have
property $P$.
\end{defn}

\begin{lem} \label{lem:Nphi}
Let $P = \{P_f\}$ be a $\Ct^k$-coordinate independent property with
$k \geq 1$.  Let $N$ be an iterative method computing $P$ with rate
$K>1$.  Fix an $f\colon \rn \rightarrow \reals$ and $\xstar \in P_f$.
Let $\phi$ be a $\Ct^k$-diffeomorphism of open
subsets of $\rn$ whose image contains $\xstar$.
Then the iteration function
\begin{equation} \label{eq:tNf}
\bar{N}_f = \phi \circ N_{f \circ \phi} \circ \phi^{-1}
\end{equation}
converges locally with rate $K$ to $\xstar$.
\end{lem}
\begin{proof}
Since $\phi^{-1}(\xstar) \in P_{f \circ \phi}$, $N_{f \circ \phi}$
converges locally to $\phi^{-1}(\xstar)$ with rate $K$.  Furthermore,
$\phi$ is bi-Lipschitz about $\phi^{-1}(\xstar)$ in a suitably small
neighbourhood of $\phi^{-1}(\xstar)$ because $k\geq 1$.  \reflem{ecnv}
completes the proof.
\end{proof}

\reflem{Nphi} suggests the coordinate-adapted viewpoint used to
generalise the Newton method to manifolds may prove beneficial in
more general contexts.  The remainder of this section elicits this
idea.

The convergence proofs for the coordinate adapted Newton method in
\refsec{cad} and the generalised Newton method in \refsec{nmm}
relied on essentially just two properties of the Newton iterate:
invariance to affine coordinate changes, and a lower bound on the radius
of convergence of the Newton iterate $N_{f \circ \phi}$ in terms
of the second-order behaviour of $\phi$.  Here, radius of
convergence refers to $\bar\rho$ in \refsec{rate_c}.

Affine invariance was
exploited partially for convenience --- it meant only parametrisations
$\phi_p$ with $\phi_p(0_p)=p$ and $D\phi_p(0_p)=I$ were needed
--- and partially to allow the Newton iterate \refeq{newt}
to be applied unambiguously to the abstract vector space
$T_pM$.  The reason for using tangent spaces was again for
convenience. It made it easier to exploit smoothness when possible.
\refsec{gtc} discussed this in detail.

Affine invariance was also used for re-centring.  In \refeq{Ef},
the change of coordinate transformations $\phi_x$ do not change the
point $x$, that is, $\phi_x(x)=x$.  This is perhaps the most natural
choice for $\phi_x(x)$ as it does not shift the space unnecessarily.
When lifting an iteration function to a manifold, such a choice is
no longer possible.  The proposed solution was to choose a distinguished
point of $\rn$, the origin, and always apply the Newton iteration
function at this distinguished point; see \refeq{gni}.  The
invariance of the Newton method to shifts made this inconsequential.

If the iterative method $N_f$ is not shift-invariant then re-centring
it at each iteration may alter its behaviour.  It is therefore
necessary to study the re-centred iterate $x_{k+1} = \theta_{x_k}
\circ N_{f \circ \theta_{x_k}}(0)$ where $\theta_x(y) = x+y$.
Equivalently, $N$ can be replaced by its re-centred version
$\tilde N_f(x) = \theta_x \circ N_{f \circ \theta_x}(0)$ which
is shift-invariant: $\theta_z \circ \tilde N_{f \circ \theta_z} \circ
\theta_z^{-1}(x) = \tilde N_f(x)$.  \textit{Henceforth the iterative
method $N_f$ is assumed to be shift-invariant.}

If the iterative method $N_f$ is affine-invariant then smooth
parametrisations $\phi\colon TM \rightarrow M$ can be used to lift $N_f$
to manifolds in the same way the Newton method was lifted.  Otherwise,
parametrisations from $M \times \rn$ rather than $TM$ need be
considered if a global approach is taken.  A simpler and more general
alternative is to construct parametrisations locally, as
in \refsec{whim}, leading to
path-dependent lifts of $N_f$.

Henceforth, a local viewpoint is adopted because determining how
to make parametrisations constructed locally fit together globally
is a topological problem unrelated to local convergence properties
and which needs no addressing if a path-dependent lift is adequate.

In one sense, lifting a shift-invariant iterative method to a
manifold locally about a point is straightforward.

\begin{prop}
\label{prop:lift}
As in \refdef{Pbar}, let $M$ be a $\Ct^k$-smooth manifold, let
$(U,\varphi)$ be a chart on $M$, and let $\bar P$ be the lift to
$M$ of a $\Ct^{k'}$-coordinate independent property $P$, with $k'
\leq k$.  Let $N$ an iterative method of order at most $k$ that
computes $P$ with rate $K$.  Define $E_f(p) = \varphi^{-1} \circ
\theta_{\varphi(p)} \circ N_{f \circ \varphi^{-1} \circ
\theta_{\varphi(p)}}(0)$ for $p \in U$, where $\theta_x(y)=x+y$.
Then $E_f$ computes $\bar P$ on $U$ with rate $K$, meaning for any
cost function $f\colon M \rightarrow \reals$ and any $\pstar \in U \cap
\bar P_f$, the iteration function $E_f$ converges locally with rate
$K$ to $\pstar$ (\refdef{conv}).
\end{prop}
\begin{proof}
Let $(V,\psi)$ be such that $\psi(\pstar) \in P_{f \circ \psi^{-1}}$.
Shrink $U$ if necessary so that $\pstar \in U \subset V$. Since $P$
is coordinate invariant, the diffeomorphism $\psi \circ
\varphi^{-1}$ can be used to show $\varphi(\pstar) \in P_{f \circ
\varphi^{-1}}$.  By \refdef{conv}, it suffices to study 
$\tilde N(x) = \varphi \circ E_f \circ \varphi^{-1}(x) =
\theta_x \circ N_{(f \circ \varphi^{-1}) \circ \theta_x}(0) =
\theta_x \circ N_{(f \circ \varphi^{-1}) \circ \theta_x} \circ
\theta_x^{-1}(x) = N_{(f \circ \varphi^{-1})}(x)$ where the last
equality follows from shift-invariance.  Since $\xstar = \varphi(\pstar)
\in P_{f \circ \varphi^{-1}}$ and $N$ computes $P$, $\tilde N$
converges locally to $\xstar$ with rate $K$, as required.
\end{proof}

Allowing more flexibility than afforded by \refprop{lift} is desirable
for two reasons: \refsec{cad_disc} explained how customised
parametrisations can improve performance for certain classes of
cost functions, and \refsec{recent} gave techniques for adapting
parametrisations to geometric features of the manifold.

The most general way found for lifting a Newton method to a manifold
is \refeq{gni}.  Furthermore, the use of $\psi$ in \refeq{gni} is
an add-on: if $E_f(p) = \phi_p \circ N_{f \circ \phi_p}(0_p)$
converges then \refeq{gni} will also converge with the same rate
provided $\psi$ is a sufficiently good approximation to $\phi$.
All that remains then is to understand when $E_f(p) = \phi_p \circ
N_{f \circ \phi_p}(0)$ computes $\bar P$ with rate $K$ given that
$N$ computes $P$ with rate $K$, as in \refprop{lift}.
Note that here, $\phi$ is defined on $U \times \rn$ where $U$ is
an open subset of $M$, and $\phi_p(0)=p$.
It is also necessary for $\phi_p$ to be a local diffeomorphism
about 0, that is, a genuine change of coordinates.

As in the proof of \refprop{lift}, consider $\tilde N(x) =
\varphi \circ E_f \circ \varphi^{-1}(x)$. 
Using the fact that
$(\varphi \circ \phi_{\varphi^{-1}(x)} \circ \theta_{-x})(x)=x$ for
all $x$, this becomes
\begin{equation}
\label{eq:tNx}
\tilde N(x) = \sigma_x
    \circ N_{\bar f \circ \sigma_x}
    \circ \sigma_x^{-1}(x),\quad
    \sigma_x = \varphi \circ \phi_{\varphi^{-1}(x)} \circ \theta_{-x},\quad
    \bar f = f \circ \varphi^{-1}.
\end{equation}
Assume $N_{\bar f}$ converges with rate $K$ to $\xstar$.
Then $\tilde N(\xstar) = \xstar$.  The trick for seeing how $\tilde
N$ converges locally to $\xstar$ is to use \refeq{tNf} to remove
the coordinate change $\sigma_\xstar$ from \refeq{tNx}.  Precisely,
the iterative method $N_g$ in \refeq{tNx} is replaced by the iterative
method $\bar N = \sigma_\xstar \circ N_{g \circ \sigma_\xstar}
\circ \sigma_\xstar^{-1}$.  By \reflem{Nphi}, this change will not alter
the rate $K$ of convergence provided $K > 1$.  (Recall from
\refapp{conv} that the $K=1$ case is more delicate.)  Thus,
\begin{equation}
\label{eq:final}
\tilde N(x) = \psi_x \circ \bar N_{\bar f \circ \psi_x}
    \circ \psi_x^{-1}(x),\quad
    \psi_x = \varphi \circ \phi_{\varphi^{-1}(x)} \circ
    \theta_{\xstar-x} \circ \phi_{\varphi^{-1}(\xstar)}^{-1} \circ \varphi^{-1}.
\end{equation}
As arranged, $\psi_\xstar$ is the identity.  If $N$, and hence $\bar
N$, is reasonably nice then the radius of convergence --- equivalently,
the constant $\kappa$ in \refsec{rate_c} --- associated with
\refeq{final} should remain bounded if $\psi_x$ remains sufficiently
close to the identity.  Indeed, all \refeq{final} is doing is
applying $\bar N$ to the cost function $\bar f$ in the coordinate
system determined by $\psi_x$.  At the end of the day, lifting
iterative methods to manifolds relies on this one simple principle:
that the iterative method be robust to changes of coordinates.


\section{Conclusion} \label{sec:concl}

The Newton method \refeq{newt} is traditionally lifted to manifolds
by endowing the manifold with a Riemannian structure and using
\refeq{riem_n}.  This strategy provides limited insight and may
have a high computational cost when implemented.  This motivates
the study, from first principles, of lifting iterative methods from
Euclidean space to manifolds.

Coordinate changes play a central role.  Changing coordinates at
each iteration is a novel yet easily understood and applied technique
for enhancing the performance of iterative methods in Euclidean
space (\refsec{cad_disc}).  Robustness to coordinate changes is key
to lifting iterative methods to manifolds in useful ways (\refsec{it}).

Newton methods on manifolds are defined customarily as iteration
functions $E_f\colon M \rightarrow M$.  This is unnecessarily restrictive;
allowing $E_f$ to depend on past history leads to path-dependent
Newton methods (\refsec{whim}), and a change in focus from devising
parametrisations to devising transformations (\refsec{recent}).
The simplifications stemming from this generalisation are a consequence
of eliminating the need for local lifts to agree globally; 
global agreement is a topological problem with little bearing on
the computational problem of iteratively finding a critical point.
While smooth global lifts of the Newton method always exist, global
agreement may not be possible if the iterative method is not invariant
to affine transformations and the manifold is non-parallelisable
(\refsec{gtc}).

The proposed framework for lifting the Newton method to manifolds
is perhaps the most general one possible (\refsec{nmm}): the condition
on the cost function in \refth{nmm} agrees with the necessary and
sufficient condition in \refth{Neuclid} for the Euclidean case, and
it is difficult to see how conditions C1--C2 can be weakened (see
\refsec{canewt}).


\appendix

\section{Rate of Convergence of Iterates on Manifolds}
\label{app:conv}

Prior to this work\footnote{The main results of this paper were
obtained in 2004--2005 and communicated privately to colleagues.},
it was natural to define convergence with respect to a Riemannian
metric.  The belief that Newton methods should not depend on any
Riemannian geometry led to the following.  Compared
with~\cite[Section 4.5]{Absil:2008uy}, the lemmata here are careful
to ensure the iterates do not fall outside the domain of definition
of the iteration function.

Convergence rates are not preserved by arbitrary homeomorphisms.
A sufficient condition for rates $K>1$ is the following.

\begin{lem}
\label{lem:ecnv}
Let $N$ be an iteration function on $\rn$ which converges locally to $\xstar$ with
rate $K>1$ and constant $\kappa$.  Let $U$ be a neighbourhood of $\xstar$
and $\phi\colon U \rightarrow V \subset \rn$ a bi-Lipschitz homeomorphism
about $\xstar$, meaning there exist positive constants $\alpha,\beta
\in \reals$ such that
\begin{equation}
\label{eq:bi}
\forall x \in U,\quad
\frac1\alpha\,\| x - \xstar \| \leq \| \phi(x) - \phi(\xstar) \| \leq
\beta\,\| x - \xstar \|.
\end{equation}
Then $\tilde{N} = \phi \circ N \circ \phi^{-1}$ converges locally
to $\phi(\xstar)$ with rate $K$ and constant $\alpha^K\beta\kappa$.
\end{lem}
\begin{proof}
As noted in \refsec{rate_c}, since $N$ converges locally to $\xstar$, for all
sufficiently small balls $B$ centred at $\xstar$, $N$ is
defined on $B$, and $x \in B$ implies $N(x) \in B$ and
$\| N(x) - \xstar \| \leq \kappa\,\|x-\xstar\|^K$.  Choose such
a $B$ contained in $U$.
Since $\phi$ is a homeomorphism, $Y = \phi(B)$ is a non-empty open
subset of $V$.  If $y \in Y$ then $\tilde{N}(y)$ is well-defined
and contained in $Y$, and
$\| \tilde{N}(y) - \phi(\xstar) \|
\leq \beta\,\| N(\phi^{-1}(y)) - \xstar \|
\leq \beta \kappa\, \| \phi^{-1}(y) - \xstar \|^K
\leq \alpha^K \beta \kappa\, \| y - \phi(\xstar) \|^K.$
\end{proof}

A significantly stronger condition is required if $K=1$.
One such example is the following.

\begin{lem}
\label{lem:lin}
Let $N$ be an iteration function on $\rn$ converging locally to
$\xstar$ at a linear rate.  Let $U$ be a neighbourhood of $\xstar$
and $\phi\colon U \rightarrow V \subset \rn$ a $\Ct^1$-diffeomorphism
whose differential $D\phi$ at $\xstar$ is proportional to the
identity.  Then $\tilde{N} = \phi \circ N \circ \phi^{-1}$ converges
locally to $\phi(\xstar)$ at a linear rate.
\end{lem}
\begin{proof}
Let $\gamma \in \reals$ be such that $D\phi(\xstar)\cdot\xi = \gamma
\xi$ for $\xi \in \rn$.  Note $\gamma \neq 0$ because $\phi$ is a
diffeomorphism.  Since $\phi(x) - \phi(\xstar) = \gamma(x-\xstar) +
r(x)$ where $\lim_{x \rightarrow \xstar} \| r(x) \| / \| x - \xstar \|
= 0$, by shrinking $U$ to become a sufficiently small neighbourhood
of $\xstar$, it can be arranged for \refeq{bi} to hold with
$\beta = |\gamma| + \epsilon$ and $\alpha = \frac1{|\gamma|-\epsilon}$
for any $\epsilon > 0$.  The result follows from \reflem{ecnv} by
choosing $\epsilon$ so that $\alpha \beta \kappa < 1$, where $\kappa <
1$ is the constant associated with $N$.
\end{proof}

The above suggests the following definition.  An iteration function
$E\colon M \rightarrow M$ on an $n$-dimensional manifold $M$ is said to
converge locally with rate $K \geq 1$ to $\pstar$ with respect to
the homeomorphism $\varphi\colon W \subset M \rightarrow V \subset \rn$,
where $\pstar \in W$, if $\varphi \circ E \circ \varphi^{-1}$, as
an iteration function on $\rn$, converges locally with rate $K$ to
$\varphi(\pstar)$.

If $K=1$ or $M$ is only a topological manifold, there is no
distinguished choice of homeomorphism $\varphi$ with respect to
which convergence can be defined.  If $M$ is a $\Ct^1$-manifold,
$K>1$ and an iterate converges with respect to one coordinate chart
$\varphi$ then \reflem{ecnv} implies it converges with respect to
any other chart $\psi$.  (Proof: If $N = \varphi \circ E \circ
\varphi^{-1}$ converges then, since $\psi \circ \varphi^{-1}$ is
$\Ct^1$ and hence bi-Lipschitz on a possibly smaller domain, $\psi
\circ E \circ \psi^{-1} = (\psi \circ \varphi^{-1}) \circ (\varphi
\circ E \circ \varphi^{-1}) \circ (\psi \circ \varphi^{-1})^{-1}$
converges too.) \refdef{conv} affords a coordinate
independent definition of rate of convergence.

\begin{defn}
\label{def:conv}
An iteration function $E\colon M \rightarrow M$ on a $\Ct^1$-differentiable
manifold converges locally with rate $K > 1$ to $\pstar \in M$ if
there exists a coordinate chart $\varphi\colon W \rightarrow V \subset
\rn$ defined on a neighbourhood of $\pstar$ such that $\varphi \circ
E \circ \varphi^{-1}$ converges locally with rate $K$ to $\varphi(\pstar)$
as an iteration function on $\rn$.
\end{defn}


\section{Local Parametrisations}
\label{app:cond}

The normalisation $\phi_x(x)=x$ used in \refsec{canewt} does not
generalise well to the manifold setting. \refsec{nmm} implicitly
introduced $h_x(y) = \phi_x(x+y)$, thereby changing the normalisation
to $h_x(0) = x$.  Properties H1 to H3 of \refsec{nmm} are the analogues
of properties P1 and P2 in \refsec{canewt}.

Choosing $h_x(y) = x + y + y^2$ if $x$ is rational and $h_x(y) = x
+ y - y^2$ if $x$ is irrational exemplifies H1--H3 do not imply
continuity of $h$.  Conversely, $h$ being $\Ct^1$-smooth and
satisfying H1 and H2 need not imply H3.

\begin{exmp}
Let $\alpha\colon \reals \rightarrow \reals$ be a $\Ct^2$-smooth (or
even $\Ct^\infty$-smooth) bump function satisfying: $0 \leq \alpha(t)
\leq 1$; $\alpha(t) = \alpha'(t) = 0$ for $t \not\in (1/2,1)$;
$\alpha(3/4) = 1$.  Let $h(x,y) = x + y + x^{-1/2} \alpha(y/x^2) y^2$
if $x > 0$ and $h(x,y)=x+y$ otherwise.  Then differentiation shows
that $h(x,y)$ is $\Ct^1$-smooth in $(x,y)$.  Furthermore, $h_x(0)=x$,
$Dh_x(0) = 1$ and $D^2h_x(0) = 0$.  Therefore, H1 and H2 are
satisfied, but H3 is not; if $x_n \rightarrow 0$ with $x_n > 0$ and
$y_n = (3/4)x_n^2$ then $(h_{x_n}(y_n) - x_n - y_n)y_n^{-2} \rightarrow
\infty$.
\end{exmp}

Nevertheless, a corollary of \reflem{P} is that $h$ being $\Ct^2$-smooth,
or even just $D^2h_x(y)$ being continuous in $(x,y)$, suffices for
H1 to imply H2 and H3.

\begin{lem}
\label{lem:P}
If, for $x,y \in B(0;\rho)$, $D^2h_x(y)$ is
bounded in $(x,y)$ and continuous in $y$ (that is, for each $x$,
$h_x(y)$ is $\Ct^2$-smooth in $y$) then $h$ satisfying \Hi{\rho} implies
it satisfies \Hii{\rho} and \Hiii{\rho}.
\end{lem}
\begin{proof}
Let $\alpha = \sup_{x,y \in B(0;\rho)} \| D^2h_x(y)\|$; then
\Hii{\rho} is satisfied.  Taylor's theorem implies $h_x(y) = x + y +
\frac12 D^2h_x(x+t(y-x)) \cdot (y-x,y-x)$ for some $t \in [0,1]$.
Thus, \Hiii{\rho} holds with $\beta = \alpha/2$.
\end{proof}

If $h(y) = y + t^3\sin(1/t)$ then $|h(y)-y| \leq |y|^2$ whenever $|y|\leq1$,
however, $D^2h(0)$ does not exist.  This puts \reflem{PP} into
context.

\begin{lem}
\label{lem:PP}
If $h$ satisfies \Hiii{\rho} then it satisfies \Hi{\rho}, and if
additionally $D^2h_x(0)$ exists for $x \in B(0;\rho)$ then $h$
satisfies \Hii{\rho} (with $\alpha=2\beta$).
\end{lem}
\begin{proof}
That \Hiii{\rho} implies \Hi{\rho} is clear.
If $D^2h_x(0)$ exists, it is known that
\begin{equation}
\lim_{\| y \| \rightarrow 0}
\| h_x(y) - 2 h_x(0) + h_x(-y) - D^2h_x(0)\cdot (y,y) \|\,\|y\|^{-2} = 0.
\end{equation}
Thus, for any $\epsilon > 0$ there is a $\delta > 0$ such that $\|
h_x(y)-2x+h_x(-y)-D^2h_x(0)\cdot(y,y)\| \leq \epsilon\|y\|^2$
whenever $\|y\| \leq \delta$.  Then $\|D^2h_x(0)\cdot(y,y)\| \leq
\epsilon\|y\|^2 + \| h_x(y)-x-y \| + \|h_x(-y)-x-(-y)\| \leq
(\epsilon + 2\beta)\|y\|^2$, proving the result; both sides
scale as $\|y\|^2$ and $\epsilon>0$ was arbitrary.
\end{proof}

\reflem{comp} asserts that \Hiii{\rho} is preserved under second-order
changes to $h$; the straightforward proof is omitted.

\begin{lem}
\label{lem:comp}
For some $\rho > 0$, assume $h$ satisfies \Hiii{\rho}.  If there
exists a $\gamma \in \reals$ such that $\tilde h$ satisfies $\|
h_x(y) - \tilde h_x(y) \| \leq \gamma \|y\|^2$ whenever $x,y \in
B(0;\rho)$ then $\tilde h$ satisfies \Hiii{\rho}.
\end{lem}

The following two technical lemmata will be required in subsequent
proofs; \reflem{lip} is well-known.

\begin{lem}
\label{lem:lip}
Given $g\colon \reals^n \rightarrow \reals^m$ and $\delta > 0$, define
$L = \sup_{z \in B(0;\delta)} \|Dg(z)\|$ and $M = \sup_{z \in
B(0;\delta)} \frac12 \|D^2g(z)\|$.  If $g$ is $\Ct^1$-smooth on
$B(0;\delta)$ and $L$ is finite then $\|g(x)-g(y) \| \leq L \|x-y\|$,
and if $g$ is $\Ct^2$-smooth on $B(0;\delta)$ and $M$ is finite
then $\|g(x)-g(y)-Dg(y)\cdot (x-y)\| \leq M\|x-y\|^2$
and $\| Dg(x) - Dg(y) \| \leq 2M\|x-y\|$ for $x,y \in
B(0;\delta)$.  If $g$ is $\Ct^1$-smooth on
$\overline{B}(0;\delta)$, meaning it is $\Ct^1$-smooth on 
an open set $U \supset
\overline{B}(0;\delta)$, then $L$ is finite, and $M$ is finite if
$g$ is $\Ct^2$-smooth on $\overline{B}(0;\delta)$.
\end{lem}

\begin{lem}
\label{lem:hiii}
Fix a dimension $n$.  Given scalars $\rho_1, \rho_2, \beta_1, G,
L, M > 0$, there exist $\rho, \beta > 0$ such that, for any 
$\bar h\colon B_n(0;\rho_1) \rightarrow \reals^n$ satisfying $\|\bar h(y) - y \|
\leq \beta_1 \|y\|^2$ for $y \in B(0;\rho_1)$, and for any $g\colon
B_n(0;\rho_2) \rightarrow \reals^n$ that is a $\Ct^2$-diffeomorphism
onto its image and satisfies $g(0)=0$, $\| [Dg(0)]^{-1}
\| \leq G$, $\| Dg(0) \| \leq L$ and $\sup_{z \in B(0;\rho_2)}
\frac12 \| D^2g(z) \| \leq M$, it follows that $\tilde
h(y) = (g \circ \bar h \circ [Dg(0)]^{-1})(y)$ is defined for $y \in
B(0;\rho)$ and satisfies $\| \tilde h(y) - y \| \leq \beta \| y \|^2$.
\end{lem}
\begin{proof}
For brevity, define $A = [Dg(0)]^{-1}$.
By successively shrinking $\rho > 0$ as required, the following
requirements can be met for all $y \in B(0;\rho)$:
$\| A y \| \leq \rho G < \rho_1$;
$\| \bar h(A y) - A y \| \leq \beta_1 \|A y\|^2 \leq \beta_1 G^2 \|y\|^2$;
$\| A^{-1} \bar h(A y) - y \| \leq \beta_1 L G^2 \|y\|^2$;
$\| \bar h(A y) \| \leq (1+\beta_1 \rho G) G\|y\| < \rho_2$;
$\| g(\bar h(A y)) - A^{-1} \bar h(A y) \| \leq M \| \bar h(A y) \|^2
    \leq M G^2 (1+\rho\beta_1 G)^2 \|y\|^2$ (\reflem{lip});
and finally $\| \tilde h(y) - y \| \leq \| \tilde h(y) - A^{-1} \bar h(A y) \|
    + \| A^{-1} \bar h(A y) - y \| \leq M G^2 (1+\rho\beta_1 G)^2 \|y\|^2
    + \beta_1 L G^2 \|y\|^2$.
Importantly, an appropriate value of $\rho$ can be determined as a
function of the other scalars and does not depend on $g$
or $\bar h$.  Similarly, $\beta = M G^2 (1+\rho\beta_1 G)^2 + \beta_1 L
G^2$ suffices.
\end{proof}

In certain situations, such as in \refsec{sym}, $h_x$ is constructed
from transformed versions of a prototype $\bar h$, as in \reflem{symh}.

\begin{lem}
\label{lem:symh}
Let $\bar h\colon \reals^n \rightarrow \reals^n$ restricted to $B(0;\rho_1)$
satisfy $\| \bar h(y) - y \| \leq \beta \|y\|^2$ for some $\beta
\in \reals$. Assume $D^2\bar h(0)$ exists. Define $h_x(y) = g_x
\circ \bar h ([Dg_x(0)]^{-1} \cdot y)$ where, for each $x \in B(0;\rho_3)
\subset \reals^n$, $g_x\colon \reals^n \rightarrow \reals^n$ restricted
to $B(0;\rho_2)$ is a $\Ct^2$-diffeomorphism satisfying $g_x(0)=x$,
$\| Dg_x(0) \| \leq L$ and $\| [Dg_x(0)]^{-1} \| \leq G$ for some
$G, L \in \reals$.  Assume $M = \sup_{x \in B(0;\rho_3),\,y
\in B(0;\rho_2)} \frac12 \|D^2g_x(y)\|^2 < \infty$. Here,
$\rho_1,\rho_2,\rho_3 > 0$.  Then $h$ satisfies H1, H2 and H3.
\end{lem}
\begin{proof}
Since $D^2 \bar h(0)$ exists and $g_x$ is $\Ct^2$-smooth, 
$D^2 h_x(0)$ exists. By \reflem{PP}, it suffices to
prove $h$ satisfies H3.  Fix $x \in B(0;\rho_3)$. Choose $\rho$ and
$\beta$ as in \reflem{hiii} with $g(z) = g_x(z) - x$. Then $\|
h_x(y) - x - y \| = \| (g \circ \bar h \circ [Dg(0)]^{-1})(y) - y
\| \leq \beta \|y\|^2$ whenever $\|y\| < \rho$.  Therefore $h$
satisfies \Hiii{\min\{\rho,\rho_2\}}.
\end{proof}

Properties H1--H3 are preserved under a change of coordinates.

\begin{lem}
\label{lem:coch}
Let $g\colon \reals^n \rightarrow \reals^n$ restricted to $B(0;\rho_3)$
be a $\Ct^2$-diffeomorphism onto its image, with $g(0)=0$.  Given
a function $h\colon \reals^n \times \reals^n \rightarrow \reals^n$,
define $\tilde h_x(y) = g \circ h_{g^{-1}(x)}( D(g^{-1})(x) \cdot
y )$.  If $h$ satisfies H1 and H2 then $\tilde h$ satisfies H1 and
H2. If $h$ satisfies H3 then $\tilde h$ satisfies H3.
\end{lem}
\begin{proof}
Assume first that $h$ satisfies \Hiii{\rho_1}.  Choose a $\rho_2$
such that $0 < \rho_2 < \frac{\rho_3}2$ and $B(0;\rho_2) \subset
g(B(0;\min\{\rho_1,\frac{\rho_3}2\}))$.  Fix an $x \in B(0;\rho_2)$.
Define $\bar h(y) = h_{g^{-1}(x)}(y) - g^{-1}(x)$ and $\tilde g(z)
= g(z+g^{-1}(x))-x$.  Note that $\bar h(y)$ is well-defined for
$\|y\| < \rho_1$ and $\tilde g(z)$ is well-defined for $\|z \| <
\rho_2$.  Then \reflem{hiii} is applicable, with $\tilde g$ replacing
$g$.  (By shrinking $\rho_3$ if necessary, it can be assumed the
derivatives of $\tilde g$ are uniformly bounded.) In particular, there
exist $\rho$ and $\beta$, independent of $x$, such that $\|\tilde
h_x(y)-x-y\| = \| (\tilde g \circ \bar h \circ [Dg(0)]^{-1})(y) -
y \| \leq \beta \|y\|^2$ whenever $y \in B(0;\rho)$.  Therefore
$\tilde h$ satisfies \Hiii{\min\{\rho,\rho_2\}}, as required.

Next, assume $h$ satisfies H1 and H2 (but not necessarily H3).  It
is reasonably clear that $\tilde h_x(y)$ has a sufficiently large
domain of definition required for $\tilde h_x(0)$, $D\tilde h_x(0)$
and $D^2\tilde h_x(0)$ to exist in a neighbourhood of $x=0$.  Explicit
calculations, using the chain rule to compute derivatives, verify
that $\tilde h$ satisfies H1 and H2.
\end{proof}

It is remarked that the tedious nature of the last few proofs
comes from the necessity of ensuring the transformed $h$ has a valid
domain of definition.  This is a consequence of the standing
assumption that $h$ itself need not be defined on the whole of
$\reals^n \times \reals^n$.  This becomes important when coordinate
charts on manifolds enter the picture.


\section{Further Results on the Generalised Newton Method} 
\label{app:fr}

\subsection{Intrinsic Conditions}

Condition \refeq{nmmf} does not depend on
the choice of coordinates.  

\begin{prop}
\label{prop:intHess}
In \refth{nmm}, if \refeq{nmmf} holds, it holds with respect to any
$\Ct^2$-chart $(\tilde U, \tilde\varphi)$ with $\tilde\varphi(\pstar)=0$
and $\tilde U$ sufficiently small.
\end{prop}
\begin{proof}
Referring to \refth{nmm}, let $(\tilde U,\tilde\varphi)$ be a chart
with $\tilde\varphi(\pstar)=0$ and choose $\rho>0$ so that $h =
\varphi \circ \tilde\varphi^{-1}$ is well-defined on
$\overline{B}(0;\rho)$.  Then $\Hess_{f \circ \tilde\varphi^{-1}}(x)
= \Hess_{\widehat f \circ h}(x) = A_x^T \Hess_{\widehat f}(h(x))
A_x + G_x$ where $A_x$ and $G_x$ are the matrix representations of
$Dh$ and $(Df \circ h) D^2h$ respectively.  Since $Dh$ and $Df \circ
h$ are $\Ct^1$-smooth and $D^2h$ is continuous, there exist constants
$\alpha,\beta$ such that $\|A_x-A_0\|\leq \alpha\|x\|$ and $\|G_x\|
\leq \beta\|x\|$ whenever $x \in \overline{B}(0;\rho)$.  Similarly,
from \refeq{nmmf} and Taylor series arguments, there exists a
constant $\gamma$ such that $\| [ \Hess_{\widehat f}(h(x)) -
\Hess_{\widehat f}(0) ] A_0 x \| \leq \| [ \Hess_{\widehat f}(h(x))
- \Hess_{\widehat f}(0) ] h(x) \| + \| [ \Hess_{\widehat f}(h(x))
- \Hess_{\widehat f}(0) ] (h(x) - A_0 x) \| \leq \gamma \|x\|^2$
whenever $x \in \overline{B}(0;\rho)$.  Shrink $\tilde U$ to equal
$\tilde\varphi^{-1}(B(0;\rho))$.  The result follows by noting
\begin{multline}
\| [\Hess_{f \circ \tilde\varphi^{-1}}(x) 
    - \Hess_{f \circ \tilde\varphi^{-1}}(0)] x \| 
\leq \| [A_x^T \Hess_{\widehat f}(h(x)) A_x  - A_x^T \Hess_{\widehat
    f}(h(x)) A_0] x \| + \\
\| [A_x^T \Hess_{\widehat f}(h(x)) A_0 - A_0^T
    \Hess_{\widehat f}(h(x)) A_0] x \| + 
\| [A_0^T \Hess_{\widehat f}(h(x)) A_0 - A_0^T
    \Hess_{\widehat f}(0) A_0] x \| +\\
\| G_x x \|.
\end{multline}
\end{proof}

Conditions C1--C2 are also intrinsic; the choice of coordinate
charts is immaterial and the conditions are preserved under diffeomorphisms.
Let $h_\star\colon TM \rightarrow TN$ denote the push-forward of tangent
vectors induced by a map $h\colon M \rightarrow N$ between manifolds;
$h_\star(v_p) = Dh(p)\cdot v_p$.

\begin{prop}
\label{prop:Cinv}
Let $\phi,\psi\colon TM \rightarrow M$ satisfy C1--C2.  Then 
about any point $p \in M$, C1 and C2 hold with respect to any $\Ct^2$-chart
$(\tilde U, \tilde\varphi)$ with $\tilde\varphi(p)=0$.
Furthermore, if $h\colon M \rightarrow N$ is a $\Ct^2$-diffeomorphism
of manifolds then the induced maps $\tilde\phi = h \circ \phi
\circ h_\star^{-1}$ and $\tilde\psi = h \circ \psi \circ
h_\star^{-1}$ satisfy C1--C2.
\end{prop}
\begin{proof}
Let $h\colon M \rightarrow N$ be a $\Ct^2$-diffeomorphism.  Fix $p \in
M$.  Let $(\tilde U, \tilde\varphi)$ be a $\Ct^2$-chart on N with
$\tilde\varphi \circ h(p) = 0$.  It will be shown $\widehat{\tilde\phi}
= \tilde\varphi \circ \tilde\phi \circ \tau_{\tilde\varphi}^{-1}$
satisfies H1 and H2, and $\widehat{\tilde\psi} = \tilde\varphi \circ
\tilde\psi \circ \tau_{\tilde\varphi}^{-1}$ satisfies H3.
This proves the second part of the lemma. The first part then follows
by letting $h\colon M \rightarrow M$ be the identity map.

Let $(U,\varphi)$ be a $\Ct^2$-chart on $M$ with $\varphi(p)=0$ and
such that $\widehat\phi=\varphi \circ \phi \circ \tau_\varphi^{-1}$
satisfies H1 and H2, and $\widehat\psi=\varphi \circ \psi \circ
\tau_\varphi^{-1}$ satisfies H3.  Let $g = \tilde\varphi \circ h
\circ \varphi^{-1}$; it is a $\Ct^2$-diffeomorphism from $\varphi(U
\cap h^{-1}(\tilde U))$ to $\tilde\varphi(h(U) \cap \tilde U)$ and
$\widehat{\tilde\psi}(x,y) = g \circ \widehat\psi_{g^{-1}(x)} \circ
D(g^{-1})(x) \cdot y$.  Apply \reflem{coch} to conclude
$\widehat{\tilde\psi}$ satisfies H3.  Analogously, \reflem{coch}
implies $\widehat{\tilde\phi}$ satisfies H1 and H2.
\end{proof}

\subsection{Sufficient Conditions}
\label{app:sc}

Conditions C1--C2 are readily satisfied by $\Ct^2$-smooth
parametrisations.  In this case $M$ must be $\Ct^3$-smooth.
If $M$ were only $\Ct^2$-smooth then $\phi\colon TM \rightarrow M$ at
best can be $\Ct^1$-smooth because $TM$ is only a $\Ct^1$-manifold.

\begin{lem}
\label{lem:Mthree}
Let $M$ be a $\Ct^3$-manifold.  If $\phi$ is $\Ct^2$-smooth and,
for all $p \in M$, $\phi_p(0_p) = p$ and $D\phi_p(0_p) = I$, then
C1 holds.  If $\psi$ is $\Ct^2$-smooth and, for all $p \in M$,
$\psi_p(0_p) = p$ and $D\psi_p(0_p) = I$, then C2 holds.
\end{lem}
\begin{proof}
Follows from \reflem{P}.
\end{proof}

\begin{rem}
\label{rem:Mthree}
Since C1--C2 are local in nature (\reflem{rp}), it suffices in \reflem{Mthree}
for $\phi$ and $\psi$ to be smooth on a neighbourhood of the zero
section of $TM$.
\end{rem}

Conditions C1--C2 are preserved under restriction to submanifolds.

\begin{lem}
\label{lem:embed}
Let $i\colon N \rightarrow M$ be a $\Ct^2$-embedding of $N$ in $M$, with
$i_\star\colon TN \rightarrow TM$ the induced push-forward of tangent
vectors.  Let $\phi,\psi\colon TM \rightarrow M$ be parametrisations of
$M$ satisfying C1--C2, and $\tilde\phi,\tilde\psi\colon TN \rightarrow
N$ parametrisations of $N$ satisfying $\phi \circ i_\star = i \circ
\tilde\phi$ and $\psi \circ i_\star = i \circ \tilde\psi$. Then
$\tilde\phi,\tilde\psi$ satisfy C1--C2.
\end{lem}
\begin{proof}
From \refprop{Cinv}, it suffices to assume $N \subset M$. Then
$\tilde\phi$ and $\tilde\psi$ are simply the restrictions of $\phi$
and $\psi$ to $TN$. The result follows by observing that if $h\colon
\reals^n \times \reals^n \rightarrow \reals^n$ in the definitions
of H1--H3 is restricted to $\reals^m \subset \reals^n$ then H1--H3
would continue to hold.
\end{proof}

One way to express precisely the local nature of C1--C2 
is with the aid of a Riemannian metric on $M$.

\begin{lem}
\label{lem:rp}
Let $\phi,\psi\colon TM \rightarrow M$ satisfy C1--C2 where $M$ is a
$\Ct^2$-Riemannian manifold.  Let $r\colon M \rightarrow (0,\infty)$ be
a possibly discontinuous function.  Assume $\tilde\phi,\tilde\psi\colon
TM \rightarrow M$ satisfy $\tilde\phi(v_p) = \phi(v_p)$ and
$\tilde\psi(v_p) = \psi(v_p)$ whenever $\| v_p \| < r(p)$.  Then
$\tilde\phi$ satisfies C1.  If $\inf_{p \in K} r(p) > 0$ for any
compact $K \subset M$ then $\tilde\psi$ satisfies C2.
\end{lem}
\begin{proof}
Fix $p \in M$ and let $\varphi$ and $\rho$ be such that C1 and C2
hold.  Then $B(0;\rho) \times B(0;\rho)$ is in the image of
$\tau_\varphi$; let $V$ be its pre-image.  For any $\bar r > 0$,
the set $V_{\bar r} = \{ v_p \in V \mid \| v_p \| < \bar r\}$ is
open, hence $\tau_\varphi(V_{\bar r})$ is open too.

Choose an $x \in B(0;\delta)$ and let $\bar r = r(\varphi^{-1}(x))$.
There exists a $\delta_x > 0$ such that $(x,B(0;\delta_x))
\subset \tau_\varphi(V_{\bar r})$.  Restricted to
$(x,B(0;\delta_x))$, $\varphi \circ \phi \circ \tau_\varphi^{-1}$
and $\varphi \circ \tilde\phi \circ \tau_\varphi^{-1}$ are equal.
It follows that $\phi$ satisfies C1.

Let $K = \varphi^{-1}(\overline{B}(0;\rho / 2))$ and $\bar r =
\inf_{p \in K} r(p)$.  Let $\bar\rho \in (0, \rho / 2)$ be such
that $B(0;\bar\rho) \times B(0;\bar\rho) \subset \tau_\varphi(V_{\bar
r})$.  Restricted to $B(0;\bar\rho) \times B(0;\bar\rho)$, $\varphi
\circ \psi \circ \tau_\varphi^{-1}$ and $\varphi \circ \tilde\psi
\circ \tau_\varphi^{-1}$ are equal.  It follows that $\psi$ satisfies
C2.
\end{proof}

\subsection{Embedded Submanifolds of Euclidean Space}
\label{sec:embed}

For manifolds embedded in Euclidean space, C1--C2 can be expressed
in extrinsic coordinates.

Treating $\reals^m$ as a manifold,
a parametrisation $\phi\colon T\reals^m \rightarrow \reals^m$ can be
specified by its representation $\widehat\phi\colon \reals^m \times
\reals^m \rightarrow \reals^m$ with respect to the identity chart,
denoted $\phi = \widehat\phi \circ \tau_I$.
Given a $\Ct^2$-embedding $i\colon M \rightarrow \reals^m$, let $V_xM$
for $x \in i(M)$ denote the realisation of $T_{i^{-1}(x)}M$ as a
subspace of $\reals^m$, that is, $(x,V_xM) = \tau_I \circ
i_\star(T_{i^{-1}(x)}M)$ where $i_\star\colon TM \rightarrow T\reals^m$
is the push-forward of $i$.  (The elements of $V_xM$ are the vectors
$\gamma'(0)$ where $\gamma\colon (-\epsilon, \epsilon) \rightarrow
\reals^m$, $\gamma(0)=x$, is a curve whose image is contained
in $i(M)$.)

If $\widehat\phi(x,y)$ belongs to $i(M)$ whenever $x \in i(M)$ and
$y \in V_xM$ then it induces a parametrisation $\tilde\phi\colon TM
\rightarrow M$ given by $\tilde\phi = i^{-1} \circ \widehat\phi
\circ \tau_I \circ i_\star$.  In essence, $\widehat\phi$ maps a
point $x+y$ on the affine tangent space of $i(M)$ at $x$, to the
point $\widehat\phi(x,y)$ on $i(M)$.  This is how parametrisations
were specified in~\cite{Manton:opt_mfold}.

\begin{lem}
\label{lem:affine_embed}
Let $i\colon M \rightarrow \reals^m$ be a $\Ct^2$-embedding of a manifold
$M$. With notation as above, assume $\widehat\phi,\widehat\psi\colon
\reals^m \times \reals^m \rightarrow \reals^m$ satisfy:
$\forall z \in i(M)$,
$\exists \alpha, \beta, \rho \in \reals$ with $\rho > 0$,
$\forall x \in B(z;\rho) \cap i(M)$,
$\forall y \in B(0;\rho) \cap V_xM$,
$\widehat\phi_x(y), \widehat\psi_x(y) \in i(M)$,
$\widehat\phi_x(0)=x$,
$D\widehat\phi_x(0) \cdot y = y$,
$\| D^2 \widehat\phi_x(0) \cdot (y,y) \| \leq \alpha \|y\|^2$,
$\| \widehat\psi_x(y) - x - y \| \leq \beta \|y\|^2$.
Then the parametrisations $\tilde\phi,\tilde\psi\colon TM \rightarrow M$
defined by
$\tilde\phi = i^{-1} \circ \widehat\phi \circ \tau_I \circ i_\star$ and
$\tilde\psi = i^{-1} \circ \widehat\psi \circ \tau_I \circ i_\star$
satisfy C1 and C2 of \refsec{nmm}.
\end{lem}
\begin{proof}
For $x \in i(M)$, let $P_x\colon \reals^m \rightarrow V_xM$ denote
Euclidean projection onto $V_xM$.  Extend $\widehat\phi$ by defining
$\widehat\phi(x,y)=x+y$ for $x \not\in i(M)$, and $\widehat\phi(x,y)=
\widehat\phi(x,P_x(y)) + y - P_x(y)$ for $x \in i(M)$ and $y \not\in
V_xM$.  Extend $\widehat\psi$ similarly.  Then $\phi=\widehat\phi
\circ \tau_I$ and $\psi=\widehat\psi \circ \tau_I$ satisfy C1--C2.
(Fix $p \in \reals^m$ and define $\varphi(x) = x-p$.  Note $\tau_I
\circ \tau_\varphi^{-1}(x,y) = (x+p,y)$.  Hence $\varphi \circ \phi
\circ \tau_\varphi^{-1}(x,y) = \widehat\phi(x+p,y)-p$.  Same for
$\psi$.  It is readily verified the assumptions in the proposition
ensure H1, H2 and H3 are satisfied.) Hence, from \reflem{embed},
$\tilde\phi$ and $\tilde\psi$ satisfy C1--C2.
\end{proof}


\section*{Acknowledgements}

This work was funded in part by the Australian Research Council.
Special thanks to Dr Jochen Trumpf for insightful and thought-provoking
discussions during the preliminary stages of this paper, and to the
two anonymous reviewers for excellent guidance on improving the
presentation.


\bibliographystyle{abbrv}
\bibliography{algmconv}

\begin{thebibliography}{10}

\bibitem{Absil:2008uy}
P.~A. Absil, R.~Mahony, and R.~Sepulchre.
\newblock {\em {Optimization algorithms on matrix manifolds}}.
\newblock Princeton University Press, Princeton, NJ, 2008.

\bibitem{Absil:2002js}
P.~A. Absil, R.~Mahony, R.~Sepulchre, and P.~Van~Dooren.
\newblock {A Grassmann-Rayleigh quotient iteration for computing invariant
  subspaces}.
\newblock {\em SIAM Review. A Publication of the Society for Industrial and
  Applied Mathematics}, 44(1):57--73, 2002.

\bibitem{Absil:2012eb}
P.~A. Absil and J.~Malick.
\newblock {Projection-Like Retractions on Matrix Manifolds}.
\newblock {\em Siam Journal on Optimization}, 22(1):135--158, 2012.

\bibitem{Absil:2004bd}
P.~A. Absil, R.~Sepulchre, P.~Van~Dooren, and R.~Mahony.
\newblock {Cubically convergent iterations for invariant subspace computation}.
\newblock {\em SIAM Journal on Matrix Analysis and Applications}, 26(1):70--96,
  2004.

\bibitem{Adler:2002cc}
R.~L. Adler, J.-P. Dedieu, J.~Y. Margulies, M.~Martens, and M.~Shub.
\newblock {Newton's method on Riemannian manifolds and a geometric model for
  the human spine}.
\newblock {\em Ima Journal of Numerical Analysis}, 22(3):359--390, 2002.

\bibitem{Alvarez:2008hv}
F.~Alvarez, J.~Bolte, and J.~Munier.
\newblock {A unifying local convergence result for Newton's method in
  Riemannian manifolds}.
\newblock {\em Foundations of Computational Mathematics}, 8(2):197--226, 2008.

\bibitem{Argyros:2007cw}
I.~K. Argyros.
\newblock {An improved unifying convergence analysis of Newton's method in
  Riemannian manifolds}.
\newblock {\em Journal of Applied Mathematics and Computing}, 25(1-2):345--351,
  2007.

\bibitem{bk:Deuflhard:newton}
P.~Deuflhard.
\newblock {\em Newton Methods for Nonlinear Problems: Affine Invariance and
  Adaptive Algorithms}.
\newblock Springer Series in Computational Mathematics. Springer, 2004.

\bibitem{Edelman:1998ei}
A.~Edelman, T.~A. Arias, and S.~T. Smith.
\newblock {The Geometry of Algorithms with Orthogonality Constraints}.
\newblock {\em SIAM Journal on Matrix Analysis and Applications},
  20(2):303--353, Jan. 1998.

\bibitem{Ferreira:2002uj}
O.~Ferreira and B.~Svaiter.
\newblock {Kantorovich's theorem on Newton's method in Riemannian Manifolds}.
\newblock {\em Journal of Complexity}, 18(1):304--329, 2002.

\bibitem{Gabay:optm}
D.~Gabay.
\newblock Minimizing a differentiable function over a differentiable manifold.
\newblock {\em Journal of Optimization Theory and Applications},
  37(2):177--219, 1982.

\bibitem{Helmke:1994ec}
U.~Helmke and J.~B. Moore.
\newblock {\em {Optimization and dynamical systems}}.
\newblock Communications and Control Engineering Series. Springer-Verlag London
  Ltd., London, 1994.

\bibitem{bk:Hirsch:diff}
M.~W. Hirsch and S.~Smale.
\newblock {\em Differential Equations, Dynamical Systems, and Linear Algebra}.
\newblock Academic Press, 1974.

\bibitem{bk:Kantorovich:fn_analysis}
L.~Kantorovich and G.~Akhilov.
\newblock {\em Functional Analysis in Normed Spaces}.
\newblock Fizmatgiz, Moscow, 1959.

\bibitem{Manton:opt_mfold}
J.~H. Manton.
\newblock Optimisation algorithms exploiting unitary constraints.
\newblock {\em IEEE Transactions on Signal Processing}, 50(3):635--650, March
  2002.

\bibitem{Manton:OG}
J.~H. Manton.
\newblock Optimisation geometry.
\newblock In K.~H{\"u}per and J.~Trumpf, editors, {\em Mathematical System
  Theory --- Festschrift in Honor of Uwe Helmke on the Occasion of his Sixtieth
  Birthday}, pages 261--274. CreateSpace, 2013.

\bibitem{bk:Ortega:iter_soln}
J.~M. Ortega and W.~C. Rheinboldt.
\newblock {\em Iterative Solution of Nonlinear Equations in Several Variables}.
\newblock Academic Press, 1970.

\bibitem{bk:Polak:opt}
E.~Polak.
\newblock {\em Optimization: Algorithms and Consistent Approximations}.
\newblock {Springer-Verlag}, 1997.

\bibitem{Shub:1986ub}
M.~Shub.
\newblock {Some remarks on dynamical systems and numerical analysis}.
\newblock In {\em Dynamical systems and partial differential equations
  (Caracas, 1984)}, pages 69--91. Univ. Simon Bolivar, Caracas, 1986.

\end{thebibliography}



\end{document}